\documentstyle{article}
\newcommand{\hqs}{\hspace{0.25in}}

\newcommand{\hm}{\hspace{5mm}}

\newcommand{\un}{\underline}
\newcommand{\ov}{\overline}
\newcommand{\sm}{\setminus}
\newcommand{\back}{\backslash}
\newcommand{\al}{\alpha}
\newcommand{\bet}{\beta}
\newcommand{\gam}{\gamma}
\newcommand{\del}{\delta}
\newcommand{\eps}{\epsilon}
\newcommand{\kap}{\kappa}
\newcommand{\lam}{\lambda}

\newcommand{\sig}{\sigma}
\renewcommand{\th}{\theta}
\newcommand{\Th}{\Theta}

\newcommand{\Cc}{\mbox{\bf C}}
\newcommand{\Rr}{\mbox{\bf R}}
\newcommand{\bt}{\begin{tabbing}}
\newcommand{\et}{\end{tabbing}}
\newcommand{\be}{\begin{equation}}
\newcommand{\ee}{\end{equation}}
\newcommand{\ds}{\displaystyle}

\newtheorem{Lemma}{Lemma}

\newtheorem{Corollary}{Corollary}
\newtheorem{Definition}{Definition}
\newtheorem{Proposition}{Proposition}
\newtheorem{Theorem}{Theorem}

\begin{document}
\title{Convergence of regular approximations to the spectra of singular
 fourth order Sturm-Liouville problems.}
\author{Malcolm Brown \thanks{Department of Computer Science, University of 
    Wales -- Cardiff, P.O. Box 916, Cardiff CF2 3XF} 
        \\ Leon Greenberg \thanks{Dept. of Mathematics, University of Maryland, 
        College Park, Maryland MD 20742} 
        \\ Marco Marletta \thanks{Dept. of Mathematics, University of Leicester, 
        University Road, Leicester LE1 7RH}}
\date{January 1996}
\maketitle
\begin{abstract}
We prove some new results which justify the use of interval truncation as
a means of regularising a singular fourth order Sturm-Liouville problem
near a singular endpoint. Of particular interest are the results in the
so called lim-3 case, which has no analogue in second order singular
problems.
\end{abstract}

\section{Introduction}\label{section:1}
In 1978, Bailey Gordon and Shampine \cite{kn:BSG} released a code (SLEIGN)
for computing eigenvalues of Sturm-Liouville problems. This code was remarkable
both for its reliability and for the fact that it was able to handle singular
as well as regular problems, with a minimum of user input. SLEIGN's strategy
for dealing with problems having singular endpoints was to truncate the interval
near a singular endpoint, thereby regularising the problem. At the time, no 
rigorous
proofs were given for the universal validity of this strategy. Recently, 
however,
a new code (SLEIGN2) has appeared, capable of dealing with much more general
singular second order problems: in particular, SLEIGN2 can deal with 
non-Friedrichs
boundary conditions near a lim-2 singular endpoint, and can deal with cases
where there is an infinite sequence of eigenvalues tending to $-\infty$. The
development of SLEIGN2 followed the work of Bailey, Everitt, Weidmann and Zettl
\cite{kn:BEWZ} in proving rigorously the types of spectral convergence which 
could be expected from the interval truncation process near singular endpoints.

In 1994, Greenberg and Marletta produced a code (SLEUTH) for solving fourth 
order
regular Sturm-Liouville problems \cite{kn:LGMM2}. This was the culmination of
three years of work on higher order self-adjoint ODEs \cite{kn:LG1, kn:LGMM1} 
and numerical
solution of eigenproblems for Hamiltonian systems \cite{kn:MM}. The new code
gave sufficient increases in speed over the Hamiltonian systems code 
\cite{kn:MM}
for it to be feasible to solve singular fourth order problems; such problems
were treated on
an heuristic basis in \cite{kn:LGMM2}. The purpose of the present work is to
prove 
results similar to those produced by Bailey,
Everitt, Weidmann and Zettl \cite{kn:BEWZ}, but in the context of fourth order
problems. For lim-2 and lim-4 singular endpoints we use methods which are
direct adaptations to the fourth order case of the methods of \cite{kn:BEWZ}.
However, for the lim-4 case, a new complication is the existence of complex as
well as real boundary conditions.
For the lim-3 case new difficulties are present, as this case does not arise
for second order problems and cannot be treated by any of the methods described
in \cite{kn:BEWZ}. We overcome these difficulties by using the oscillation 
theory
described in \cite{kn:LGMM1}.

\section{ODE theory for the fourth order problem}\label{section:background}
A fourth order Sturm-Liouville equation is an equation of the form

\be \ell y = \lambda y, \;\;\; a < x < b,\label{eq:mmintro1} \ee
where $\ell $ is a differential operator of the form
\be \ell y := \frac{1}{w(x)}\left\{ [(p(x)y'')'-(s(x)y']'+q(x)y\right\}.
 \label{eq:mmintro2} \ee
Here $1/p$, $s$, $q$ and $w$ are locally $L^1$ in $(a,b)$, with $p$ and
$w$ positive almost everywhere, $s$ and $q$ real-valued. The endpoints $a$ and 
$b$ may be finite
or infinite. The endpoint $x=a$ is {\bf regular} if it is finite and if there 
exists
$a' >a$ such that $1/p$, $s$, $q$ and $w$ are in $L^1(a,a')$; a similar 
definition
holds for $x=b$. Any other sort of endpoint is called {\bf singular}.

Let $L^2(a,b;w)$ be the space of functions $f(x)$ on $(a,b)$ such that
\[  \int_a^b |f(x)|^2w(x)\,dx < \infty . \]
$L^2(a,b;w)$ is a Hilbert space with inner product
\[ \langle f,g\rangle =\int_a^bf(x)\overline{g(x)}w(x)\,dx, \]
and norm $\|f\| = \sqrt{\langle f,f\rangle}$.
\noindent We shall say that a function $f$ is
{\bf square integrable at $a$} if, for some $\beta \in (a,b)$, 
\[ \int_a^{\beta}|f(x)|^2w(x)\,dx < \infty. \]
Square integrability at $b$ is defined similarly.

Given a function
$y$ we define the {\bf quasi-derivatives} 
\be y^{[0]} = y, \;\; y^{[1]} = y', \;\; y^{[2]}
  = py'', \;\; y^{[3]} = -(py'')'+sy'. \label{eq:mmqd}  \ee
\noindent  These quasi-derivatives were introduced for scalar $2n^{th}$ order
problems by Naimark \cite{kn:N}; see also Everitt and Zettl \cite{kn:EZ}
and Zettl \cite{kn:Zettl} for further information on quasi-differential
operators.
\begin{Definition} \label{def:maxdom}
The {\bf maximal domain} $D_{max}$ is the set of functions $y$ whose
quasi-derivatives $y^{[0]}$, $y^{[1]}$, $y^{[2]}$ and $y^{[3]}$ are all 
absolutely
continuous, and for which $y$ and $\ell y$ lie in $L^2(a,b;w)$.
The {\bf maximal operator} $L_{max}$ is the operator defined by
$L_{max}y = \ell y$ on the domain $D(L_{max}) = D_{max}$.
\end{Definition}

\begin{Definition}
The {\bf pre-minimal domain} $C_{min}$ is the set of all functions in $D_{max}$
having compact support in $(a,b)$. The {\bf minimal domain} $D_{min}$ is the
closure of $C_{min}$ in the graph norm
\[ \| y \|_{G}^2 = \langle y,y \rangle + \langle \ell y, \ell y \rangle. \]
The minimal operator $L_{min}$ is the restriction of $L_{max}$ to $D_{min}$.
\end{Definition}

It is not difficult to see that $L_{min}$ is the adjoint of $L_{max}$. Every
self-adjoint extension of $L_{min}$ is a restriction of $L_{max}$
to some domain $D$ between $D_{min}$ and $D_{max}$.
It may happen that $L_{max} = L_{min}$, in which case both
are self-adjoint and $\ell$ has only this one self-adjoint realisation. However,
$L_{max}$ is not generally self-adjoint. 
It is known that a self-adjoint extension $L$ of $L_{min}$ has a domain
$D(L) \subset D_{max}$ which is determined by certain boundary conditions
at the endpoints.
In order to know how many boundary conditions are required, we need
the following classification of endpoints.
\begin{Definition}\label{mmdef5}
An endpoint is said to be of {\bf lim-$p$ type} if the space of solutions
of the differential equation $\ell y = iy $ which are square integrable at
that endpoint has dimension $p$.
\end{Definition}

For fourth order Sturm Liouville problems, the only possibilities are
$p=2$, $p=3$, and $p=4$. The lim-2 case is analogous to the limit-point
case for second order equations; the lim-4 case is analogous to limit-circle.
Lim-3 has no second order analogue. Regular endpoints are of lim-4 type.
These endpoint types should not be confused with the deficiency indices
of the minimal operator: for example, if both endpoints are of lim-2
type then the deficiency indices are zero.

Boundary conditions are imposed by means of the Lagrangian form, which 
we now define.
\begin{Definition}\label{mmdef6} The {\bf Lagrangian form}
$[f,g]$ of two functions $f$ and $g$ in $D_{max}$ is defined
by

\[ [f,g](x):= 
\left\{f^{[0]}(x)\overline{g^{[3]}(x)}+f^{[1]}(x)\overline{g^{[2]}(x)}\right\} - 
 \left\{f^{[3]}(x)\overline{g^{[0]}(x)}+f^{[2]}(x)\overline{g^{[1]}(x)}\right\} 
. \]
with the notation of (\ref{eq:mmqd}) for quasiderivatives.
\end{Definition}
Using another notation, 
\[ [f,g](x) = u_f^T\ov{v}_g - v_f^T\ov{u}_g,   \]
where
\[ u_f = \left(\begin{array}{l} f^{[0]} \\ f^{[1]} \end{array}\right), \;\;\;
   v_f = \left(\begin{array}{l} f^{[3]} \\ f^{[2]} \end{array}\right), \]
and the superscript $T$ indicates the transposed matrix.  Integration by 
parts shows that for\hfill\linebreak 
\mbox{$f, g \in D_{max}$} and \mbox{$a < \al < \beta < b$},
\be \int_{\al}^{\beta}(\ell f)\ov{g}w\,dx -\int_{\al}^{\beta}f(\ov{\ell g})w\,dx
    = [f,g](\beta) - [f,g](\al).  \label{eq:lg1}  \ee
This shows that the limits
\[ [f,g](a) := \lim_{x\rightarrow a+}[f,g](x), \;\;\; 
 [f,g](b) := \lim_{x\rightarrow b-}[f,g](x) \]
both exist.

The domain $D(L)$ of a self-adjoint extension of $L_{min}$ is determined
by boundary conditions of the form $[y,\psi](b) - [y,\psi](a) = 0$, where
$\psi \in D_{max} \back D_{min}$.  We are interested in separated boundary
conditions, which are of the form $[y,\phi](a) = 0$ or $[y,\psi](b) = 0$.
For such conditions, we may assume that $\phi = 0$ in a left neighborhood
of $b$, and $\psi = 0$ in a right neighborhood of $a$.

At a lim-2 endpoint no boundary conditions are required. At
a lim-3 endpoint, one boundary condition is required. Suppose,
for example, that $x=a$ is a lim-3 endpoint. Then we can choose
any $\phi \in D_{max}\backslash D_{min}$, such that $[\phi,\phi](a) = 0$, 
and impose a condition
\be [y,\phi](a) = 0. \label{eq:mm29} \ee
All valid boundary conditions at $x=a$ would have this form.
Now suppose that $x=a$ were of lim-4 type. Then we would have
to choose two maximal domain functions $\phi_1$ and $\phi_2$ 
with $[\phi_1,\phi_1](a)=[\phi_1,\phi_2](a)=[\phi_2,\phi_2](a)=0$, 
and such that no (non-trivial) linear 
combination of $\phi_1$ and $\phi_2$ is in $D_{min}$. With these
functions we would then impose two boundary conditions
\be [y,\phi_1](a)=0, \;\;\; [y,\phi_2](a) = 0. 
 \label{eq:mm30}
\ee

Equation (\ref{eq:lg1}) proves the following result which
will be required later.
\begin{Lemma}\label{mmlem10}
If $y_1$ and $y_2$ are two solutions of $\ell y = \lambda y$ for
the same real $\lambda$ then $[y_1,y_2]$ is constant. If $y_1$
and $y_2$ are solutions for the same complex $\lambda$ then
$[y_1,\overline{y_2}]$ is constant.
\end{Lemma}

We have now seen how to construct self-adjoint realisations $L$ of
$\ell$ by restriction of $L_{max}$ to a domain $D(L) \subset D_{max}$.
As we mentioned earlier, self-adjoint realisations of $\ell$ may also
be obtained by extension of $D_{min}$. It is particularly useful to
view $D(L)$ in this way when constructing a core for $L$.

\begin{Definition} \label{mmdef1} A set $C\subset D(L)$ is said
to be a {\bf core} of $L$ if $D(L)$ is the closure of $C$ in
the graph norm.
\end{Definition}

In essence there are only six possible combinations of endpoint
singularities: lim-2--lim-2, lim-2--lim-3, lim-2--lim-4,
lim-3--lim-3, lim-3--lim-4 and lim-4--lim-4. We shall consider
each case in turn.

\begin{Lemma} (Domains and cores of self-adjoint extensions) \label{mmlemma11}
\begin{description}
\item[lim-2--lim-2] In this case $D(L)=D_{min}$. Since $D_{min}$ is
 the closure of $C_{min}$ in the graph norm, $C_{min}$ is a 
 core of $L$.
\item[lim-2--lim-3] Suppose we have a boundary condition
 $[y,\psi](b)=0$ and suppose without loss of generality that
 $\psi$ is zero on a right neighbourhood of $x=a$. Then
 $D(L) = D_{min}\oplus\mbox{Span}(\psi)$. The set 
 $C_{min}\oplus\mbox{Span}(\psi)$ is a core of $L$.
\item[lim-2--lim-4] Suppose that we have boundary conditions
 $[y,\psi_1](b)$ $=$ $0$ $=$ $[y,\psi_2](b)$ and suppose
 that $\psi_1$ and $\psi_2$ are zero on a right neighbourhood
 of $x=a$. Then $D(L)=D_{min}\oplus\mbox{Span}(\psi_1,\psi_2)$.
 The set $C_{min}\oplus\mbox{Span}(\psi_1,\psi_2)$ is a core
 of $L$.
\item[lim-3--lim-3] Suppose that we have boundary conditions
 $[y,\phi](a)=0$ and $[y,\psi](b)=0$ and suppose that
 $\psi$ is zero on a right neighbourhood of $x=a$ while
 $\phi$ is zero on a left neighbourhood of $x=b$. Then
 $D(L)=D_{min}\oplus\mbox{Span}(\phi,\psi)$. The set
 $C_{min}\oplus\mbox{Span}(\phi,\psi)$ is a core of $L$.
\item[lim-3--lim-4] Suppose that we have boundary conditions
 $[y,\phi](a)=0$ and $[y,\psi_1](b)$ $=$ $0$ $=$ $[y,\psi_2](b)$
 and suppose that $\phi$ is zero on a left neighbourhood of
 $x=b$ while $\psi_1$ and $\psi_2$ are zero on a right neighbourhood
 of $x=a$. Then $D(L)=D_{min}\oplus\mbox{Span}(\phi,\psi_1,\psi_2)$.
 The set $C_{min}\oplus\mbox{Span}(\phi,\psi_1,\psi_2)$ is a core
 of $L$.
\item[lim-4--lim-4] Suppose that we have boundary conditions
 $[y,\phi_1](a)$ $=$ $0$ $=$ $[y,\phi_2](a)$ and $[y,\psi_1](b)$ 
 $=$ $0$ $=$ $[y,\psi_2](b)$  and suppose that $\phi_1$ and 
 $\phi_2$ are zero on a left neighbourhood of
 $x=b$ while $\psi_1$ and $\psi_2$ are zero on a right neighbourhood
 of $x=a$. Then 
 $D(L)=D_{min}\oplus\mbox{Span}(\phi_1,\phi_2,\psi_1,\psi_2)$.
 The set $C_{min}\oplus\mbox{Span}(\phi_1,\phi_2\psi_1,\psi_2)$ 
 is a core of $L$.
\end{description}
\end{Lemma}

We shall finish this section by describing some
results from the theory of regular fourth order problems which we
shall require later in this paper. We start with an alternative
description of boundary conditions for a regular endpoint. 
Given a function $y$ we define a vector ${\bf z_y}$ by
\be {\bf z_y} = (y^{[0]},y^{[1]},y^{[3]},y^{[2]})^T = 
(y,y',-(py'')'+sy',py'')^T. 
\label{eq:mmintro3} \ee
Recall that we have defined the vectors $u_y$ and $v_y$ by
\be u_y = (y^{[0]},y^{[1]})^T, \;\;\; v_y = (y^{[3]},y^{[2]})^T. 
\label{eq:mmintro3a} \ee
Define matrices $J$ and $S$ by
\be J = \left( \begin{array}{rrrr} 0 & 0 & -1 &  0 \\
                                   0 & 0 &  0 & -1 \\
                                   1 & 0 &  0 &  0 \\
                                   0 & 1 &  0 &  0 \end{array}\right),
\;\;\;
    S = \left( \begin{array}{cccc} \lambda w - q & 0 & 0 & 0 \\
                                          0 & -s & 1 & 0 \\
                                          0 & 1 & 0 & 0 \\
                                          0 & 0 & 0 & 1/p \end{array}\right).
\label{eq:mmintro4} \ee
The equation (\ref{eq:mmintro1}) is equivalent to the Hamiltonian
system 
\be J{\bf z'_y} = S {\bf z_y}. \label{eq:mmintro5} \ee
At a regular endpoint, say $x=a$, one imposes self-adjoint boundary conditions
as follows. Let $A_1$ and $A_2$ be two $2\times 2$ matrices such that the
$2\times 4$ matrix $(A_1  A_2)$ is of full rank (2) and such that
$A_1 A_2^* $ is Hermitian. Then the condition
\be A_1 u_y(a) + A_2 v_y(a) = 0  \label{eq:mmintro6} \ee
is a valid self-adjoint boundary condition at $x=a$, and every 
self-adjoint boundary condition at $x=a$ has this form for suitable
$A_1$ and $A_2$.

Associated with the matrices $A_1$ and $A_2$ is a matrix solution $Y_L$
of the Hamiltonian system. $Y_L$ is a $4\times 2$ matrix such that
$JY_L' = SY_L$ and 
\be Y_L(a) = \left( \begin{array}{r} -A_2^* \\ A_1^* \end{array} \right).  
\label{eq:lgintro}   \ee
Usually $Y_L$ is partitioned into two $2\times 2$ matrices $U_L$ and
$V_L$:
\[  Y_L = \left( \begin{array}{c} U_L \\ V_L \end{array} \right). \]
The matrix $U_L$ is singular only at isolated points, and so we may
define almost everywhere the matrix
\[ W_L = V_L U_L^{-1}. \]
It may be shown that $W_L$ is Hermitian when $\lambda$ is real. Similar 
ideas hold for a right endpoint $x=b$, although it is then usual to 
denote the relevant matrices by $Y_R$, $U_R$, $V_R$ and $W_R$.

In this paper we shall use extensively the oscillation theory developed
in \cite{kn:LG1} and summarised, for the case of fourth order problems,
in \cite{kn:LGMM1}.   
The theory was developed for real boundary conditions in \cite{kn:LG1},
but it can be shown that the same proofs extend to the case of complex
boundary conditions.
The main result we shall require is the following.
\begin{Theorem}\label{Theorem:intro}
Consider a regular fourth order problem over a finite interval
$(a,b)$ with separated self-adjoint boundary conditions at $a$
and $b$. Let
\[ \delta_L(x,\lambda) = \sum_{a < t < x}\left\{\mbox{rank deficiency
 of $U_L(t,\lambda)$} \right\}, \]
\[ \delta_R(x,\lambda) = \sum_{x < t < b}\left\{\mbox{rank deficiency
 of $U_R(t,\lambda)$} \right\}. \]
Then 
\begin{itemize}
\item $\delta_L$ and $\delta_R$ are finite;
\item given any point $c\in [a,b]$ there is an integer $\sigma(c,\lambda)$
 whose value is either 0, 1 or 2, such that the number of eigenvalues 
 of the eigenproblem strictly less than $\lambda$ is 
\be N(\lambda) = \delta_L(c,\lambda) + \delta_R(c,\lambda) + \sigma(c,\lambda);
\label{eq:mmN} \ee
\item if $U_L(c,\lambda)$ and $U_R(c,\lambda)$ are nonsingular then
\[ \sigma(c,\lambda) = \nu_{\#}(W_L(c,\lambda)-W_R(c,\lambda)),  \] 
where $\nu_{\#}(W) = \mbox{the number of negative eigenvalues of $W$}$;
\item if $c=b$ and Dirichlet boundary conditions $y(b) = 0 = y'(b)$ are imposed
 then $\delta_R = \sigma = 0$ in (\ref{eq:mmN});
\item if $c=a$ and Dirichlet boundary conditions $y(a) = 0 = y'(a)$ are imposed
 then $\delta_L = \sigma = 0$ in (\ref{eq:mmN}).
\end{itemize}
\end{Theorem}
Formulae for $\sigma(c,\lambda)$ covering all the other cases are given
in \cite{kn:LGMM1}.

\section{Operator convergence and spectral inclusion}

\setcounter{equation}{0}
\setcounter{Theorem}{0}
\setcounter{Corollary}{0}
\setcounter{Lemma}{0}
\setcounter{Definition}{0}

In this section we shall review some functional
analytic concepts and results which will be required in later
sections. This will be followed by some new results on spectral
inclusion for singular fourth order operators approximated
by regular fourth order operators.
\subsection{Functional analysis}
\noindent {\bf Notation}\, We consider a self-adjoint operator $L$
on a domain $D(L)$ in a Hilbert space $H$ with inner product
$\langle \cdot,\cdot \rangle$. We also consider a sequence of
self adjoint operators $(L_j)_{j=1}^{\infty}$ on domains $D(L_j)$
in $H$. We shall denote by $\mbox{Sp}(L)$ the spectrum of $L$
and by $\mbox{Sp}(L_j)$ the spectrum of $L_j$ for each $j$.
The notation $\| \cdot \|_G$ will denote the {\bf graph norm}
associated with $L$:
\[ \| y \|_G^2 := \langle y,y \rangle + \langle Ly,Ly \rangle,
 \;\;\; y\in D(L). \]
\begin{Definition}\label{mmdef2} The sequence $(L_j)_{j=1}^{\infty}$
is {\bf strong resolvent convergent} (SRC) to $L$ if, for
$z \in \Cc \backslash \Rr$ and $f\in H$, 
\[ (L_j-zI)^{-1}f \longrightarrow (L-zI)^{-1}f \;\;\; 
 \mbox{as $j\rightarrow\infty$}.\]
\end{Definition}
A sufficient condition for strong resolvent convergence is
given by the following result.
\begin{Lemma}\label{lemma:mmlem3}
Suppose that $C$ is a core of $L$ such that, for each $f\in C$,
there is an integer $N$ such that $f\in D(L_j)$ for all $j>N$,
and
\[ \lim_{j\rightarrow\infty}L_jf = Lf. \]
Then the sequence $(L_j)_{j=1}^{\infty}$ is SRC to $L$.
\end{Lemma}
See Reed and Simon \cite[Theorem VIII.25]{kn:Reed} for a proof.
\begin{Definition}\label{mmdef3} We say that $(L_j)_{j=1}^{\infty}$
is {\bf norm resolvent convergent} (NRC) to $L$ if, for some
$z\in \Cc \backslash \Rr$, 
\[ \lim_{j\rightarrow\infty} \| (L_j-zI)^{-1} - (L-zI)^{-1}\| 
   = 0. \]
\end{Definition}
\begin{Definition}\label{mmdef4}
\hspace{1.0in}
\begin{description}
\item[(i)] The sequence $(L_j)_{j=1}^{\infty}$ is {\bf
 spectrally inclusive} for $L$ if, for every $\lambda\in 
 \mbox{Sp}(L)$, there is a sequence $(\lambda^{(j)})_{j=1}^{\infty}$
 with $\lambda^{(j)}\in \mbox{Sp}(L_j)$ such that
 $\lim_{j\rightarrow\infty}\lambda^{(j)} = \lambda$.
\item[(ii)] The sequence $(L_j)_{j=1}^{\infty}$ is {\bf
 spectrally exact} for $L$ if it is spectrally inclusive for
 $L$ and if, given a sequence $(\lambda^{(j)})_{j=1}^{\infty}$ with
 $\lambda^{(j)}\in \mbox{Sp}(L_j)$, every limit point of
 the sequence is in $\mbox{Sp}(L)$.
\end{description}
\end{Definition}

The reason for introducing the concept of SRC lies in the following 
result (see Reed and Simon \cite[Theorem VIII.24]{kn:Reed}).
\begin{Theorem} \label{Theorem:SRCSI}
Suppose that $L$ is a self-adjoint operator on a Hilbert space $H$
and let $L_j$ be a sequence of self-adjoint operators on $H$ which
is SRC to $L$. Then $L_j$ is a spectrally inclusive sequence for
$L$. 

Suppose also that we denote by $P(\zeta_1,\zeta_2)$ the spectral
projection of $L$ associated with an interval $(\zeta_1,\zeta_2)$
whose endpoints are not in the spectrum of $L$ and let 
$P_j(\zeta_1,\zeta_2)$ 
be the corresponding spectral projection for $L_j$, for each $j$.
Then for any $f\in H$,
\[ \lim_{j\rightarrow\infty} P_j(\zeta_1,\zeta_2)f = P(\zeta_1,\zeta_2)f. \]
In particular, this implies that if $\lambda$ is an isolated eigenvalue of
multiplicity $k$ for $L$ then there must be $k$ eigenvalues (counted 
according to multiplicity) of $L_j$ which converge to $\lambda$ as 
$j\rightarrow \infty$.
\end{Theorem}

There is a corresponding connection between NRC and spectrally exact
convergence, given by the following result (see Reed and Simon
\cite[Theorem VIII.23]{kn:Reed}).
\begin{Theorem}\label{Theorem:NRCSE}
Suppose that $L$ is a self-adjoint operator on a Hilbert space $H$
and let $L_j$ be a sequence of self-adjoint operators on $H$ which
is NRC to $L$. Then $L_j$ is a spectrally exact sequence for
$L$. 

Furthermore, using the notation of Theorem \ref{Theorem:SRCSI}
for the spectral projections, if $(\zeta_1,\zeta_2)$ is an interval
whose endpoints are not in the point spectrum of $L$ then
\[ \lim_{j\rightarrow\infty} \| P_j(\zeta_1,\zeta_2) - 
 P(\zeta_1,\zeta_2) \| = 0. \]
\end{Theorem}
\vspace{2mm}

\noindent{\bf \un{Remark}}.\, It may happen that $L$ is 
unbounded below but that each $L_j$ is bounded below. In this
case it will be possible to index the eigenvalues of each
$L_j$ as, say, 
\[ \lambda_0^{(j)} \leq \lambda_1^{(j)} \leq \cdots, \]
assuming that these all exist. In this case, Theorem
\ref{Theorem:NRCSE} would tell us that for each fixed $k$,
$\lim_{j\rightarrow\infty}\lambda_k^{(j)} = -\infty$.

\subsection{Fourth order Sturm-Liouville operators: 
 interval truncation and spectral inclusion}
\label{subsubsection:sisect}
We construct regular approximations to our singular Sturm-Liouville
problems by interval truncation. We shall require two technical
lemmas to allow us to set up approximating boundary conditions near
lim-3 and lim-4 endpoints.
\begin{Lemma} \label{lemma:lgbcs} 
Suppose that the endpoint $x=a$ is of lim-3 type for the
differential expression $\ell$.  Let $\phi(x)$ be a function in
$D_{max} \setminus D_{min}$ such that $[\phi,\phi](a) = 0$.
Then there is a real function $z(x)$ such that the boundary condition
$[y,\phi](a) = 0$ is equivalent to $[y,z](a) = 0$.  A similar result
holds for the endpoint $x=b$.
\end{Lemma}

\noindent {\bf Proof}:\, We may suppose that $\phi(x) = 0$ in a 
left-neighborhood of $x=b$.  Let $\phi(x) = u(x) + iv(x)$, where $u(x)$
and $v(x)$ are real.  Then
\[ [\phi,\phi](a) = [u,u](a) + [v,v](a) + 2i[v,u](a) = 0.  \]
Therefore $[v,u](a) = 0$.  Furthermore,
\[  [u,u](a) = 0 = [v,v](a),  \]
since $u$ and $v$ are real.  This implies that
\[  [u,\phi](a) = [u,u](a) - i[u,v](a) = 0,  \]
\[  [v,\phi](a) = [v,u](a) - i[v,v](a) = 0.  \]
Thus $u$ and $v$ satisfy the boundary condition $[y,\phi](a) = 0$,
and therefore
\[ u = y_1 + c_1\phi, \ \ \ v = y_2 + c_2\phi,  \]
where $y_1, y_2 \in D_{min}$ and $c_1, c_2$ are constants.
If $c_1 = c_2 = 0$, then both $u \ {\rm and} \  v \in D_{min}$, and so
$\phi \in D_{min}$.  Therefore one of the constants, say $c_1 \neq 0$.
This implies that $u \in D_{max} \setminus D_{min}$ and the condition
$[y,\phi](a) = 0$ is equivalent to $[y,u](a) = 0$. \hfill $\Box$

\begin{Lemma} \label{lemma:mmbcs}
 Suppose that the endpoint $x=a$ is of lim-4 type
 for the differential expression $\ell$ and suppose that 
 boundary conditions $[y,\phi_1](a)$ $=$ $0$ $=$ $[y,\phi_2](a)$
 are imposed at $x=a$. Then for each real $\lambda$ there are
 solutions of the equation $\ell y = \lambda y$, say $z_1$
 and $z_2$, such that the boundary conditions at $x=a$ are
 equivalent to the conditions $[y,z_1](a)$ $=$ $0$ $=$ 
 $[y,z_2](a)$. Similar results hold for the endpoint $x=b$.
\end{Lemma}
\noindent {\bf Proof}:\, As we are interested in separated
boundary conditions we do not want to be concerned with
a possible singular endpoint at $x=b$; to this end we assume
that $x=b$ is regular, with two regular self-adjoint boundary
conditions imposed there. If this is not the case we can
always arrange for it to be so, by interval truncation, without
making any changes to the boundary conditions at $x=a$.
Consider the symmetric closed operator $L_a$  whose domain $D_a$
is the set of functions $y \in D_{max}$ such that
\begin{itemize}
\item $y$ satisfies the regular boundary conditions at $b$;
\item $[y,g](a) = 0$ for any $g$ in $D_{max}$.
\end{itemize}
Note that $L_a$ is a symmetric operator.  For if $y_1, y_2\in D_a$, 
then
\[ \int_a^b (\ell y_1)\ov{y}_2 w\,dx - \int_a^b y_1(\ov{\ell y_2}) w\,dx
    = [y_1,y_2](b) = 0, \]
since $y_1$ and $y_2$ satisfy self-adjoint boundary conditions at $x = b$.
The adjoint $L_a^*$ of $L_a$ is the maximal operator associated with
$L_a$ and its domain is the set of functions $y \in D_{max}$ such that
$y$ satisfies the regular boundary conditions at $b$.
\vskip 5pt

\noindent The deficiency index of $L_a$ is the dimension of $N(L_a^*-iI)$.
It is clearly two, because the fact that $x=a$ is of lim-4 type
means that the space of solutions which are square integrable
at $x=a$ is of dimension four; $N(L_a^*-iI)$ is the two-dimensional
subspace of these satisfying the two boundary conditions at $x=b$.
We now fix any real $\lambda$ and any point $c$ in $(a,b)$ and
construct four solutions $y_1,\ldots,y_4$ spanning the whole
solution space for $\ell y = \lambda y$, say by using the
initial conditions
\[ \left( \begin{array}{c} y_j(c) \\ y_j'(c) \\ -(py_j'')'(c)+sy_j'(c)
 \\ py_j''(c) \end{array}\right) = \mbox{\bf e}_j, 
 \;\;\; j=1,\ldots,4,\]
where the $\mbox{\bf e}_j$ are the standard unit vectors in
$\Rr^4$. For each $i$ and $j$ the Lagrangian form $[y_i,y_j](x)$
is constant, so it is easy to see that the matrix $W = ([y_i,y_j])$
is nonsingular. We now form four more functions $\tilde{y}_1,
\ldots,\tilde{y}_4$ which are functions in $D_{max}$ with
$\tilde{y}_j(x) = y_j(x)$ in a right neighbourhood of $x=a$ (which includes
the point $x=c$), and
$\tilde{y}_j(x) \equiv 0$ in a left neighbourhood of $x=b$.
The matrix of Lagrangian forms for these functions, $\tilde{W}$,
has the property that
\[ \tilde{W}(b) - \tilde{W}(a) = -W(a) = -W(c), \]
and $W(c)$ is nonsingular. By Lemma 10.2.17 of Hutson and Pym
\cite{kn:Hutson} this means that $\tilde{y}_1,\ldots,\tilde{y}_4$ are
linearly independent relative to  $D_a$, and by
Theorem 10.2.18 of \cite{kn:Hutson} any two self-adjoint boundary
conditions at $x=a$ are equivalent to boundary conditions of the form
\be [y,z_1](a) = 0, \;\;\; [y,z_2](a) = 0, \label{eq:mmlemma22} \ee
where $z_1,z_2 \in \mbox{Span}(\tilde{y}_1,\ldots,\tilde{y}_4)$.
In particular, our boundary conditions $[y,\phi_1](a)=0$
and $[y,\phi_2](a)=0$ can be replaced by conditions of the
form (\ref{eq:mmlemma22}). Finally since $\tilde{y}_j = y_j$ near
$a$ for $j=1,\ldots,4$, the proof is complete. \hfill $\Box$

\noindent{\bf \un{Remark}}.\, The functions $z_1, z_2$ in the previous lemma 
need not
be real.  There may be no real functions satisfying the lemma.

We are now ready to construct regular truncated interval
approximations to our singular problems. We approximate the
interval $(a,b)$ by intervals $(a_j,b_j)$, where 
$a_j\rightarrow a+$ and $b_j\rightarrow b-$ as $j\rightarrow\infty$.
If $x=a$ is regular then we can take $a_j=a$ for all $j$,
and if $x=b$ is regular we can take $b_j=b$ for all $b$.
Otherwise we must assume that $a_j\in (a,b)$ and $b_j\in (a,b)$
with $a_j<b_j$ for all $j$. At a regular endpoint the boundary
conditions are inherited from the original problem. The
list below describes what happens at singular endpoints.
In each case, we take for a core of $L$ the core described
in Lemma \ref{mmlemma11}.

\begin{description}
\item[lim-2--lim-2] In this case we impose boundary
 conditions $y(a_j)=0=y'(a_j)$ and $y(b_j)=0=y'(b_j)$.
 Any function in the core of $L$ will have compact support 
 in $(a,b)$ and will therefore satisfy these boundary conditions for
 all sufficiently large $j$.
\item[lim-2--lim-3] Suppose we have a boundary condition
 $[y,\psi](b)=0$ where $\psi$ is zero on a right neighbourhood
 of $x=a$.  By Lemma \ref{lemma:lgbcs} we may suppose that $\psi$
 is a real function, which implies that $[\psi,\psi](b_j) = 0$.
 Three of our boundary conditions are clear:
 we impose $y(a_j)=0=y'(a_j)$ and $[y,\psi](b_j)=0$.
 For the fourth boundary condition we impose any self
 adjoint condition at $x=b_j$. Any function in the core
 of $L$ will certainly satisfy the conditions at $x=a_j$
 for all sufficiently large $j$; also, since any function
 $y$ from the core is a multiple of $\psi$ near $b$ it will
 satisfy $[y,\psi](b_j)=0$. Finally since the remaining
 boundary condition is a valid self-adjoint boundary
 condition, it is necessarily satisfied by $\psi$,
 hence by multiples of $\psi$, hence by $y$ for all
 sufficiently large $j$.
\item[lim-2--lim-4] Suppose that we have boundary conditions
 $[y,\psi_1](b)$ $=$ $0$ $=$ $[y,\psi_2](b)$ and suppose
 that $\psi_1$ and $\psi_2$ are zero on a right neighbourhood
 of $x=a$. The boundary conditions at $x=a_j$ are as in
 the first two cases. At $x=b_j$ we would like to impose
 the conditions $[y,\psi_1](b_j)=0$ and $[y,\psi_2](b_j)=0$.
 We know that $[\psi_1,\psi_2](b)=0$ but we do not know
 that $[\psi_1,\psi_2](b_j)=0$, so unfortunately the obvious
 boundary conditions may not be self-adjoint. To get round
 this we invoke Lemma \ref{lemma:mmbcs}, which tells us
 that we can assume that $\psi_1$ and $\psi_2$ are solutions
 of the differential equation with real $\lambda$ 
 near $b$, so that $[\psi_1,\psi_1]$, $[\psi_1,\psi_2]$
 and $[\psi_2,\psi_2]$ are constant (and therefore $0$)
 near $b$ by Lemma \ref{mmlem10}.  Thus our boundary conditions 
 $[y,\psi_1](b_j)=0$ and $[y,\psi_2](b_j)=0$
 are self-adjoint. If a function $y$ is in the core of $L$
 then near $x=b$ it is a linear combination of $\psi_1$
 and $\psi_2$ and therefore satisfies the boundary conditions.
\item[lim-3--lim-3] Suppose that we have boundary conditions
 $[y,\phi](a)=0$ and $[y,\psi](b)=0$ and suppose that
 $\psi$ is zero on a right neighbourhood of $x=a$ while
 $\phi$ is zero on a left neighbourhood of $x=b$.  We may
 assume that $\phi$ and $\psi$ are real.  Then
 we impose boundary conditions $[y,\phi](a_j)=0$ and 
 $[y,\psi](b_j)=0$, plus one other self-adjoint condition
 at each of $a_j$, $b_j$. Any function in the core of $L$
 is a multiple of $\phi$ near $a$ and a multiple of $\psi$
 near $b$ and therefore satisfies the boundary conditions
 at $a_j$ and $b_j$ for all sufficiently large $j$.
\item[lim-3--lim-4] Suppose that we have boundary conditions
 $[y,\phi](a)=0$ and $[y,\psi_1](b)$ $=$ $0$ $=$ $[y,\psi_2](b)$
 and suppose that $\phi$ is zero on a left neighbourhood of
 $x=b$ while $\psi_1$ and $\psi_2$ are zero on a right neighbourhood
 of $x=a$.  We may assume that $\phi$ is real. 
 We impose the condition $[y,\phi](a_j)=0$ together
 with one other self-adjoint condition at $a_j$. For the
 conditions at $b_j$ we assume once more that $\psi_1$
 and $\psi_2$ are solutions of the differential equation with real
 $\lambda$ near 
 $b$, and we impose the conditions $[y,\psi_1](b_j)$ $=$ $0$ $=$ 
 $[y,\psi_2](b_j)$. Once more we can show that any function 
 in the core satisfies the boundary conditions for all sufficiently
 large $j$.
\item[lim-4--lim-4] Suppose that we have boundary conditions
 $[y,\phi_1](a)$ $=$ $0$ $=$ $[y,\phi_2](a)$ and $[y,\psi_1](b)$ 
 $=$ $0$ $=$ $[y,\psi_2](b)$  and suppose that $\phi_1$ and 
 $\phi_2$ are zero on a left neighbourhood of
 $x=b$ while $\psi_1$ and $\psi_2$ are zero on a right neighbourhood
 of $x=a$. We can assume also that $\phi_1$ and $\phi_2$ are
 solutions of the differential equation with real $\lambda$ near $x=a$ 
 and that $\psi_1$ and $\psi_2$ are solutions near $x=b$, and then impose
 the conditions $[y,\phi_1](a_j)$ $=$ $0$ $=$ $[y,\phi_2](a_j)$ and
 $[y,\psi_1](b_j)$ $=$ $0$ $=$ $[y,\psi_2](b_j)$ at $a_j$ and
 $b_j$ respectively. Any function in the core of $L$ will satisfy
 these conditions for all sufficiently large $j$.
\end{description}

Of course, the problem of maintaining self-adjointness of the
boundary conditions when truncating at a lim-4 endpoint do not
arise if the lim-4 endpoint happens to be regular, because then
no truncation is required.

Associated with each truncated interval problem is a self-adjoint
operator $L_j$ whose domain contains functions defined over
$[a_j,b_j]$ and satisfying appropriate boundary conditions at
$a_j$ and $b_j$. If we are to invoke the results of the previous
section on SRC to obtain spectral inclusion, then we must 
create from the operators $L_j$ new self-adjoint operators
$L'_j$ which act on functions defined over the whole
of $(a,b)$. This is an uninspiring technical process which
we shall now discuss.

We shall consider $L^2(a_j,b_j;w)$ to be a closed subspace of
$L^2(a,b;w)$ by extending the functions in $L^2(a_j,b_j;w)$
to be zero in $(a,b) \back (a_j,b_j)$.  Then we have a splitting
\[ L^2(a,b;w) = L^2(a_j,b_j;w) \oplus L^2(a_j,b_j;w)^{\perp},  \] 
where
\[  L^2(a_j,b_j;w)^{\perp} = L^2(a,a_j;w) \oplus L^2(b_j,b;w). \]
The operator $L_j$ is extended to an operator $L'_j$ on $L^2(a,b;w)$ 
by defining it to be zero on $L^2(a_j,b_j;w)^{\perp}$:
\be  L'_j = L_j \oplus \Th_j, \label{eq:mm33}   \ee
where $\Th_j$ is the zero operator on $L^2(a_j,b_j;w)^{\perp}$.
The domain of $L'_j$ is
\[ D(L'_j) = D(L_j) \oplus L^2(a_j,b_j;w)^{\perp}. \]
\noindent It is easy to check that the $L'_j$ are self-adjoint.
Unfortunately they also all possess 0 as an eigenvalue of
infinite multiplicity. This is no problem to us in the
interpretation of the spectral inclusion results so long as
we assume -- as we may -- that 0 is not an eigenvalue of $L$.

\begin{Lemma} \label{lemma:mmsi}
In each of the cases described above, the operators $L'_j$
are SRC to $L$.
\end{Lemma}
\noindent {\bf Proof}:\, We shall use Lemma \ref{lemma:mmlem3}.
Suppose that $f$ lies in the core of $L$ as described earlier.
Our construction of the $L_j$ has ensured that $f$ satisfies
the boundary conditions at $a_j$ and $b_j$ for all sufficiently
large $j$, and so $f$ lies in $D(L'_j)$. Thus we can
compute $Lf - L'_jf$, and from (\ref{eq:mm33}) we clearly
have
\[ (Lf-L'_jf)(x) = \left\{ \begin{array}{ll}
                   0 & \mbox{for $x\in (a_j,b_j)$,} \\
                   (\ell f)(x) & \mbox{for $x\in (a,b)\backslash (a_j,b_j)$}.
                   \end{array}\right. \]
Thus 
\[ \| Lf-L'_jf \|^2 = 
 \int_{a}^{a_j} w(x)|\ell f |^2(x)dx + 
  \int_{b_j}^{b} w(x)|\ell f |^2(x)dx. \]
Since $f$ lies in $D_{max}$, $\ell f$ must be square integrable
near $a$ and $b$ and so the integrals on the right hand side tend
to zero as $j\rightarrow \infty$. By Lemma \ref{lemma:mmlem3} this
completes the proof. \hfill $\Box$
\begin{Corollary}\label{Corollary:spec.inclus}
The truncated operators $L_j$ constructed above are spectrally inclusive
for $L$.
\end{Corollary}
\noindent {\bf Proof}:\, Theorem \ref{Theorem:SRCSI} implies that the $L'_j$
are spectrally inclusive for $L$, and therefore this is true for the $L_j$ also.
\hfill $\Box$
\section{The Friedrichs extension for lim-2, lim-3 and lim-4}
\label{section:Friedrichs}

\setcounter{equation}{0}
\setcounter{Theorem}{0}
\setcounter{Corollary}{0}
\setcounter{Lemma}{0}
\setcounter{Definition}{0}

We consider in this section a problem having one or two singular
endpoints. We assume that the minimal operator is bounded below,
so that it possesses a Friedrichs extension. Our objective is to
approximate the eigenvalues below the essential spectrum of the 
Friedrichs extension.

\subsection{Friedrichs boundary conditions at every singular endpoint}
Let us start by setting up the appropriate domains for the
minimal operators. Consider first the case of one singular endpoint,
say $x=b$. We assume that two regular boundary conditions are given
in the usual way at $x=a$. For our pre-minimal domain we take
 the set $C_{min}$ of all functions in the maximal domain which
satisfy the boundary conditions at $x=a$ and have compact support
in $[a,b)$. The minimal domain $D_{min}$ is the closure of $C_{min}$ 
in the graph norm. For the case of two singular endpoints, we take the
pre-minimal domain to be the set $C_{min}$ of all functions with
compact support in $(a,b)$; the minimal domain $D_{min}$ is, once
more, the closure of $C_{min}$ in the graph norm. In either case
we denote by $L_{min}$ the operator given on $D_{min}$ by
$L_{min}y = \ell y$. Since $L_{min}$ is bounded below we may assume,
by making a shift of $q(x)$ in the expression for $\ell$ if necessary,
that
\be \langle \ell y,y \rangle \geq \langle y,y \rangle \;\;\; \forall
 y \in D_{min}. \ee
We can then define the {\bf energy norm} $\| \cdot \|_{E} $ by
\be \| y \|_{E}^2 = \langle \ell y,y \rangle, \;\;\; y\in D_{min}.
 \label{eq:enorm} \ee
The energy norm is stronger than the weighted $L^2$ norm but not
as strong as the graph norm. Recall that the Friedrichs extension
is defined to be the operator whose domain is the closure of
$D_{min}$ in the energy norm (see, e.g., Dunford and Schwartz
\cite[pp. 1240-1241]{kn:Dunford}). We shall denote the Friedrichs extension by
$L_F$ and its domain by $D_F$. Obviously the energy norm defined
in (\ref{eq:enorm}) extends to $D_F$.

We now set up approximating regular problems. In the case of one
singular endpoint, which we have taken for convenience to be
$x=b$, we choose a monotone increasing sequence of points $b_j$ 
such that $a < b_j < b$ for all $j$ and $b_j\rightarrow b$ as 
$j\rightarrow\infty$. We set up regular problems on the truncated
intervals $[a,b_j]$ which inherit their regular boundary conditions
at $x=a$ from the original problem and which have Dirichlet boundary
conditions $y = y' = 0$ at $x=b_j$. For two singular endpoints we
require a second (monotone decreasing) sequence of points $a_j$
such that $a < a_j < b_j < b$
for all $j$ and $a_j\rightarrow a$ as $j\rightarrow\infty$. In this
case we also impose Dirichlet conditions $y=y'=0$ at the $a_j$.
We denote by $L_j$ the resulting regular operators and by
$\mu_{k}^{(j)}$ their eigenvalues, for $k=0,1,2,\ldots$; the 
eigenvalues of $L_F$ will be denoted by $\mu_k$, when they
exist.

Finally, as a notational convenience, we shall define sets
$H_j$ as follows. In the case of one singular endpoint at $x=b$,
$H_j$ will be the set of maximal domain functions on $[a,b)$
which satisfy the boundary conditions at $x=a$ and have
compact support in $[a,b_j)$; in the case of two singular endpoints
$H_j$ will be the set of maximal domain functions having
compact support in $(a_j,b_j)$. 

The following two lemmas are well known. The first is a standard
result from the theory of regular Sturm-Liouville problems; the
second comes from Berkowitz \cite[Theorem 2.2]{kn:Berk} and the fact that $D_F$ 
is the
completion of $C_{min}$ under the energy norm.
\begin{Lemma} \label{lemma:lf1} 
For each $k$ and $j$,
\[ \mu_{k}^{(j)} = \inf_{\mbox{\em dim}(V)=k+1,\; V\subset H_j }\left\{
             \sup_{v\in V} \frac{\| v \|_{E}}{\| v\|} \right\}. \]
\end{Lemma}

\begin{Lemma} \label{lemma:lf2} Suppose that $L_F$ possesses a $k$th
eigenvalue $\mu_k$ strictly below any essential spectrum. Then
\[ \mu_{k} = \inf_{\mbox{\em dim}(V)=k+1,\; V\subset C_{min} }\left\{
             \sup_{v\in V} \frac{\| v \|_{E}}{\| v\|} \right\}. \]
\end{Lemma}

We can now prove the following useful result.
\begin{Lemma} \label{lemma:lf3} 
For each $j$ and $k$, 
\[ \mu_{k}^{(j+1)} \leq \mu_{k}^{(j)}. \]
Moreover if $\mu_k$ exists and is below the essential spectrum of
$L_F$ then
\[ \mu_{k} \leq \mu_{k}^{(j)}  \]
for all $j$ and $k$.
\end{Lemma}

\noindent{\bf Proof}:\, Since $b_j$ is monotone increasing (and
$a_j$, when required, is monotone decreasing) it is clear that
$H_{j} \subset H_{j+1}$ for all $j$, from which the first inequality
holds from Lemma \ref{lemma:lf1}. The second inequality comes from
Lemma \ref{lemma:lf2} and the observation that $H_j \subset C_{min}$
for all $j$. \hfill $\Box$
\vspace{2mm}

\noindent We can now state and prove our main result.
\begin{Theorem}\label{Theorem:FApprox} Suppose that $L_F$ possesses
a $k$th eigenvalue $\mu_k$ strictly below any essential spectrum. Then
\[ \lim_{j\rightarrow\infty}\mu_{k}^{(j)} = \mu_{k}. \]
\end{Theorem}

\noindent{\bf Proof}:\, From Lemma \ref{lemma:lf2}, given $\epsilon>0$ 
there exists a $(k+1)$-dimensional space $V_{k}^{\epsilon} \subset C_{min}$ 
such that
\be \mu_{k} \geq \sup_{v\in V_{k}^{\epsilon}}\frac{\| v\|_{E}}{\| v\|}
                                           - \epsilon. \label{eq:lf1}
\ee
Because $V_{k}^{\epsilon}$ is finite-dimensional there exists a 
positive integer $J_{k}^{\epsilon}$ such that for all $j\geq J_{k}^{\epsilon}$,
$V_k^{\epsilon}\subset H_j$. Thus from (\ref{eq:lf1}),
\[ \mu_{k} \geq \inf_{\mbox{dim}(V)=k+1,\; V\subset H_j} 
            \left\{\sup_{v\in V}\frac{\| v\|_{E}}{\| v\|}\right\}
                                                       - \epsilon \]
for all $j\geq J_{k}^{\epsilon}$. The right hand side of this inequality
is just $\mu_{k}^{(j)}-\epsilon$, so 
\[ \mu_{k} \geq \mu_{k}^{(j)} - \epsilon \]
for all $j\geq J_{k}^{\epsilon}$. Combining this with the second inequality
in Lemma \ref{lemma:lf3} yields the required result. \hfill $\Box$

\subsection{Two singular endpoints with Friedrichs boundary conditions
 at one endpoint only}
Theorem \ref{Theorem:FApprox} deals with two cases: the case of one 
singular endpoint,
and the case of two singular endpoints where we want the extension
corresponding to the Friedrichs boundary conditions at both ends.
When we come to considering problems with one lim-2 endpoint
and one other singular (either lim-3 or lim-4) endpoint, we shall
find it useful to have a result on eigenvalue convergence in which
we approximate a problem with two singular ends by a problem with
one singular end. We shall therefore mention the case where a problem
has one endpoint, say $x=a$, at which we impose Friedrichs boundary
conditions (if required, i.e. if the endpoint is of lim-3 or lim-4
type), and where the remaining endpoint $x=b$ is singular, 
but where no truncation is effected there. The point $x=b$ may
be of lim-2 type, in which case no boundary conditions are required
there, or it may be of lim-3 or lim-4 type.
In the lim-3 case there will be a boundary condition at
$x=b$ of the form $[y,\psi](b)=0$; in the lim-4 case there will
be two boundary conditions at $x=b$, of the form $[y,\psi_1](b) = 0
= [y,\psi_2](b)$. We denote by $L_F$ the self-adjoint operator arising
from the equation, the Friedrichs boundary conditions at $x=a$, when
required, and any other requisite boundary conditions at $x=b$.
The domain $C_{min}$ of our preminimal operator we now define as
follows:
\begin{description}
\item[lim-2 at $b$:] $y\in C_{min}$ if and only if
 $y$ is a maximal domain function and there exists
 $\alpha\in (a,b)$ for which $y$ has support in $(\alpha,b)$.
\item[lim-3 at $b$:] $y\in C_{min}$ if and only if
 $y$ is a maximal domain function, $[y,\psi](b)=0$, and there exists
 $\alpha\in (a,b)$ for which $y$ has support in $(\alpha,b)$.
\item[lim-4 at $b$:] $y\in C_{min}$ if and only if
 $y$ is a maximal domain function, $[y,\psi_1](b)=0=[y,\psi_2](b)$,
 and there exists $\alpha\in (a,b)$ for which $y$ has support in $(\alpha,b)$.
\end{description}
The minimal domain is, as usual, the closure of $C_{min}$ in the
graph norm.

We choose a sequence of points $a_j\rightarrow a$, with $a_j\in (a,b)$ for
each $j$, and we set up truncated interval eigenproblems on the
intervals $(a_j,b)$ by imposing the regular Dirichlet boundary
conditions
\[ y(a_j) = 0 = y'(a_j). \]
We do not truncate at $b$; the boundary conditions at $b$, if 
required, remain unchanged, so our approximating Sturm-Liouville
problems are now also singular. We denote the associated operators
by $L_j$. To recover the results of the lemmas above we need
to re-define the sets $H_j$. 
\begin{description}
\item[lim-2 at $b$:] In this case we take $H_j$ to be the set of
 maximal domain functions with support in $(a_j,b)$.
\item[lim-3 at $b$:] In this case we take $H_j$ to be the set of
 maximal domain functions with support in $(a_j,b)$ such that
 $[y,\psi](b)=0$.
\item[lim-4 at $b$:] In this case we take $H_j$ to be the set of
 maximal domain functions with support in $(a_j,b)$ such that
 $[y,\psi_1](b)=0=[y,\psi_2](b)$.
\end{description}

Denote by $\mu_k^{(j)}$ the eigenvalues of $L_j$ and by $\mu_k$
the eigenvalues of $L_F$. We now have the following result.

\begin{Theorem}\label{Theorem:FApprox2}
 Suppose that the minimal operator is bounded below and suppose
that $\mu_k$ lies strictly
below any essential spectrum of $L_F$. Then
\be \lim_{j\rightarrow\infty}\mu_k^{(j)} = \mu_k. \label{eq:mm22} \ee
\end{Theorem}

Before proving this theorem we require a technical lemma.

\begin{Lemma}\label{Lemma:split}
Let ${\cal C}(\cdot )$ denote the essential spectrum of an operator.
Then ${\cal C}(L_j)\subset {\cal C}(L_F)$ for all $j$.
\end{Lemma}

\noindent {\bf Proof}:\, This is a simple application of the operator
splitting technique of Akhiezer and Glazman 
\mbox{\cite[p. 520]{kn:Akhiezer}}.
Let $D_j$ be the set of all functions $f \in D(L_F)$ with quasiderivatives
$f^{[0]}(a_j) = f^{[1]}(a_j) = f^{[2]}(a_j) = f^{[3]}(a_j) = 0$; and let
$M_j$ be the restriction of $L_F$ to $D_j$.  
Let ${\cal A}_j = D_j \cap L^2(a,a_j;w)$, 
\mbox{${\cal B}_j = D_j \cap L^2(a_j,b;w)$},
and let $A_j,\ B_j$ be the restrictions of $L_F$ to ${\cal A}_j,\ {\cal B}_j$,
respectively.  Then $D_j = {\cal A}_j \oplus {\cal B}_j$ and
$M_j = A_j \oplus B_j$.  The truncated operator $L_j$ is a self-adjoint 
extension of $B_j$.  Let $K_j$ be any self-adjoint extension of $A_j$.  
Then $K_j \oplus L_j$ and $L_F$ are two self-adjoint
extensions of $M_j$.  By \cite[Theorem 1, \S 105]{kn:Akhiezer}, 
\[ {\cal C}(L_F) = {\cal C}(K_j\oplus L_j) = {\cal C}(K_j)\cup {\cal C}(L_j). \]
Therefore ${\cal C}(L_j)\subset {\cal C}(L_F)$.\hfill $\Box$
\vskip 5pt

\noindent {\bf Proof of Theorem \ref{Theorem:FApprox2}}:\, Because 
$\mu_k$ lies below the essential spectrum of $L_F$, Lemma \ref{lemma:lf2} 
still holds. Thus, given $\epsilon>0$ 
there exists a $(k+1)$-dimensional space $V_{k}^{\epsilon}$ in the
preminimal domain $C_{min}$ such that
\be \mu_{k} \geq \sup_{v\in V_{k}^{\epsilon}}\frac{\| v\|_{E}}{\| v\|}
                                           - \epsilon. \label{eq:mm23}
\ee
Because $V_{k}^{\epsilon}$ is finite-dimensional and because each element
of $C_{min}$ belongs to all the $H_j$ for all sufficiently large $j$,
there exists a 
positive integer $J_{k}^{\epsilon}$ such that for all $j\geq J_{k}^{\epsilon}$,
$V_{k}^{\epsilon}\subset H_j$. Thus from (\ref{eq:mm23}),
\[ \mu_{k} \geq \inf_{\mbox{dim}(V)=k+1,\; V\subset H_j} 
            \left\{\sup_{v\in V}\frac{\| v\|_{E}}{\| v\|}\right\}
                                                       - \epsilon \]
for all $j\geq J_{k}^{\epsilon}$. Rearranging this result gives
\be \inf_{\mbox{dim}(V)=k+1,\; V\subset H_j} 
            \left\{\sup_{v\in V}\frac{\| v\|_{E}}{\| v\|}\right\}
 \leq \mu_k + \epsilon. \label{eq:mm24} \ee
We may assume that $\epsilon$ is so small that $\mu_k + \epsilon$ 
lies strictly below the essential spectrum of $L_F$ and hence, by Lemma
\ref{Lemma:split}, below the essential spectra of all the $L_j$; in particular 
the
left hand side of (\ref{eq:mm24}) lies strictly below the essential
spectrum of $L_j$. Standard variational theory 
(see, e.g. Berkowitz \cite[Theorem 2.2]{kn:Berk}) 
now tells us that the eigenvalue $\mu_k^{(j)}$ exists and is given by the
expression on the left hand side of (\ref{eq:mm24}), whence we
deduce that
\[ \mu_{k}^{(j)} \leq \mu_k + \epsilon \]
for all $j > J_{k}^{\epsilon}$. The inclusion $H_j \subset C_{min}$
and the now established variational expression for $\mu_k^{(j)}$ now
yield the inequality
\[ \mu_k \leq \mu_k^{(j)} \]
for all $j\geq J_{k}^{\epsilon}$. This completes the proof. \hfill $\Box$

\section{The lim-4 case}\label{section:lim4}

\setcounter{equation}{0}
\setcounter{Theorem}{0}
\setcounter{Corollary}{0}
\setcounter{Lemma}{0}
\setcounter{Definition}{0}

\subsection{Truncation at the lim-4 endpoint}

We consider here the problem where the endpoint $x=b$ is lim-4,
and the other endpoint is lim-2 or lim-3.  We shall truncate only
at the lim-4 endpoint $x=b$.  Later we shall consider double
truncations.  $L$ will denote a given self-adjoint extension
of $L_{min}$.

\begin{Lemma}\label{lemma:lgsols}
Let $x=a$ be a lim-2 or lim-3 endpoint, and $x=b$ a lim-4 endpoint.
Let $\lam \in \Cc \sm \Rr$, and let $L$ be a self-adjoint extension  
of $L_{min}$.  There are two independent solutions
$\phi_1(x), \phi_2(x)$ of the equation $\ell y = \lam y$, such that:
\begin{description}
\item[(1)]  $\phi_1, \phi_2 \in L^2(a,b;w)$;
\item[(2)]  $\phi_1$ and $\phi_2$ satisfy the boundary conditions for $L$
             at $x=a$ (if any);
\item[(3)]  $[\phi_1,\ov{\phi}_2](a) = 0$.
\end{description}
\end{Lemma}

\noindent{\bf Proof}:\, If the endpoint $x=a$ is lim-2, there exist two
independent solutions $\phi_1, \phi_2 \in L^2(a,b';w)$ for any $b' \in (a,b)$.
These solutions also lie in $L^2(b',b;w)$, since $x=b$ is a lim-4 endpoint.
Thus they lie in $L^2(a,b;w)$.  In this case there are no boundary conditions
at $x=a$, and $[\phi_1,\ov{\phi}_2](a) = 0$, since $[y_1,y_2](a) = 0$ for all
functions $y_1, y_2 \in D_{max}$.
\vskip 5pt

\noindent Suppose that $x=a$ is a lim-3 endpoint.  By Lemma \ref{lemma:lgbcs}
the boundary condition at $x=a$ is of the form $[y,z](b) = 0$, where $z(x)$ is
a real function.  This implies that if $\phi(x)$ satisfies the boundary 
condition,
then $\ov{\phi(x)}$ does also.  
\vskip 5pt

\noindent  There are three independent solutions in $L^2(a,b;w)$.  Two linearly
independent solutions $\phi_1, \phi_2$ satisfy the boundary condition at $x=a$.
Since $\ov{\phi}_2$ also satisfies the boundary condition, 
$[\phi_1,\ov{\phi}_2](a) = 0$.  \hfill $\Box$
\vskip 10pt

\begin{Lemma}\label{lemma:lg1}
Let $x=b$ be a lim-4 endpoint, and let $\lam \in \Cc$.  Suppose that
$\phi$ is a solution of $\ell y = \lam y$ such that $[y,\ov{\phi}](b) = 0$
for all solutions $y$.  Then $\phi = 0$.
\end{Lemma}

\noindent {\bf Proof}:\, For any two solutions $y_1, y_2$, \  
$[y_1,\ov{y}_2](x)$
is constant.  This defines a skew-symmetric bilinear form on the space $S$ of
solutions.  This form is nondegenerate in the sense that if $z \in S$ and
$[y,\ov{z}] = 0$ for all $y \in S$, then $z = 0$.  The nondegeneracy can be
checked as follows.  Let $a < c < b$, and for $j = 1, 2, 3, 4$, let $y_j$ be the
solution of $\ell y = \lam y$ with the following quasi-derivatives (using the
notation of (\ref{eq:mmqd})):

\[  (y_j^{[0]}(c)\ \ y_j^{[1]}(c)\ \ y_j^{[3]}(c)\ \ y_j^{[2]}(c))^T = {\bf 
e}_j,  \]

\noindent where $\{{\bf e}_1,\ {\bf e}_2,\ {\bf e}_3,\ {\bf e}_4\}$ is the 
standard basis 
for $\Cc^4$. The matrix $A = ([y_j,\ov{y}_k])$ is

\[ A = \left(\begin{array}{cc}  0 & I_2 \\ 
                             -I_2 & 0 \end{array} \right), \]
                             
\noindent where $I_2$ is the 2$\times$2 identity matrix.  Since $\det A \neq 0$,
the bilinear form $[y,\ov{z}]$ is nondegenerate.  Since $[y,\ov{\phi}] = 0$ for 
all
$y\in S$, it follows that $\phi = 0$.  \hfill $\Box$
\vskip 10pt

\begin{Definition}
Let $x=b$ be a lim-4 endpoint, and let $L$ be a self-adjoint extension
of $L_{min}$.  The boundary conditions for $L$ at $x=b$ will be 
called {\em real} if they are equivalent to conditions
$[y,z_1](b) = 0 = [y,z_2](b)$, where $z_1$ and $z_2$ are real functions.
Otherwise the boundary conditions will be called {\em complex}.
\end{Definition}

\begin{Lemma}\label{lemma:lg2}
Let $x=b$ be a lim-4 endpoint.  Let $\lam \in \Cc \sm \Rr$, and let $L$ be
a self-adjoint \mbox{extension of $L_{min}$}.  
\begin{description}
\item[(1)] There are two independent solutions $\psi_1(x),\ \psi_2(x)$  of 
the equation $\ell y = \lam y$, such that

    \begin{description}
    \item[(a)] $\psi_1,\ \psi_2 \in L^2(a',b;w)$ for any $a' \in (a,b)$;
    \item[(b)] $\psi_1$ and $\psi_2$ satisfy the boundary conditions for
               $L$ at $x=b$.
    \end{description}
    
\item[(2)] $[\psi_1,\ov{\psi}_2](b) = 0$ if and only if the boundary conditions
at $x=b$ are real.

\item[(3)] If the boundary conditions at $x=b$ are complex, then $\psi_1,\ 
\psi_2$
can be chosen so that \mbox{$[\psi_1,\ov{\psi}_2](b) = 1$}.
\end{description}
\end{Lemma}

{\bf Proof}:\, (1)  Let $S$ be the 4-dimensional space of solutions of 
$\ell y = \lam y$.  Since $x=b$ is lim-4, all functions in $S$ belong to 
$L^2(a',b;w)$ for any $a' \in (a,b)$.  There is a 2-dimensional subspace
$T \subset S$ of functions satisfying the boundary conditions
$[y,z_1](b) = 0 = [y,z_2](b)$ of $L$ at $x=b$.  Take $\psi_1,\ \psi_2$ to
be any two independent functions in $T$.
\vskip 5pt

\noindent (2)  If the boundary conditions at $x=b$ are real, then $\ov{\psi}$
satisfies these conditions if $\psi$ does.  Thus $\ov{\psi}_2$ satisfies the
boundary conditions, and so $[\psi_1,\ov{\psi}_2](b) = 0$.
\vskip 5pt

\noindent Conversely, suppose that $[\psi_1,\ov{\psi}_2](b) = 0$.  Since 
$[\psi_2,\ov{\psi}_2](b) = 0$ also, $\ov{\psi}_2$ satisfies the boundary
conditions at $x=b$.  Similarly, $\ov{\psi}_1$ satisfies these conditions.
\vskip 5pt

\noindent We claim that the boundary conditions $[y,z_1](b) = 0 = [y,z_2](b)$
for $L$ at $x=b$ are equivalent to the conditions
$[y,\psi_1](b) = 0 = [y,\psi_2](b)$. To show this, let $a<a'<b$.  For $a'<x<b$,
\begin{eqnarray} 
  \psi_1(x) & = & y_1(x) + c_1z_1(x) + c_2z_2(x), \label{eq:100}  \\
  \psi_2(x) & = & y_2(x) + d_1z_1(x) + d_2z_2(x), \label{eq:101}
\end{eqnarray}
where $y_1,\ y_2 \in D_{min}$.   By Lemma \ref{lemma:lg1} a nontrivial linear 
combination of $\psi_1$ and $\psi_2$ cannot lie in $D_{min}$.  Therefore
equations (\ref{eq:100}) and (\ref{eq:101}) can be solved for $z_1$ and
$z_2$ in terms of $\psi_1$, $\psi_2$ and functions in $D_{min}$.  This implies
that the boundary conditions are equivalent to $[y,\psi_1](b) = 0 = 
[y,\psi_2](b)$.
\vskip 5pt

\noindent Let $\psi_1 = u_1 + iv_1$, $\psi_2 = u_2 + iv_2$, where the functions
$u_k$ and $v_k$ are real.  Since $\ov{\psi}_1,\ \ov{\psi}_2$ satisfy the
boundary conditions, the real functions $u_1,\ v_1,\ u_2,\ v_2$ also satisfy
the conditions.  Two of these functions must be independent relative to 
$D_{min}$.
For otherwise $\psi_1$ and $\psi_2$ would be dependent relative to $D_{min}$.
Therefore $\psi_1$ and $\psi_2$ can be expressed as linear combinations of two
of these real functions modulo $D_{min}$, and so the boundary conditions can
be given by real functions.  \mbox{This proves (2)}.
\vskip 5pt

\noindent (3)  If the boundary conditions at $x=b$ are complex, then
$[\psi_1,\ov{\psi}_2](b) = \al \neq 0$.  The function 
$\hat{\psi}_1 = \frac{1}{\al}\psi_1$ satisfies 
$[\hat{\psi}_1,\ov{\psi}_2](b) = 1$.  Thus if $\psi_1$ is replaced by
$\hat{\psi}_1$, then the required equation is satisfied. \hfill $\Box$
\vskip 10pt

We will now calculate the Green's function for $L$.  We shall first calculate  
it for real boundary conditions at $x=b$, and then indicate the result for
complex boundary conditions.  The equation 
$\ell y - \lam y = f$ has the form:
\be [(py'')' - (sy')]' + (q-\lam w)y = fw.  \label{eq:lg40}  \ee
\noindent  This is transformed to Hamiltonian form as follows.
Corresponding to the function $y(x)$, consider the quasi-derivatives

\[ u_1 = y,\hm u_2 = y',\hm v_1 = sy'-(py'')',\hm v_2 = py'',  \]

\[u_y = \left( \begin{array}{cc}u_1\\u_2\end{array} \right),\hm
           v_y = \left( \begin{array}{cc}v_1\\v_2\end{array} \right), \hm
                  z_y = \left( \begin{array}{cc}u_y\\v_y\end{array} \right).  \]
                  
\noindent Equation (\ref{eq:lg40}) is equivalent to

\be Jz' = Sz + \hat{f}w,  \label{eq:lg41}  \ee 

\noindent where $z = z_y$, $J$ and $S$ are as in (\ref{eq:mmintro4}),
and
\[ \hat{f} = (f,0,0,0)^T. \]

\noindent Let $\phi_1, \phi_2, \psi_1, \psi_2$ be the solutions of $\ell y = 
\lam y$
from Lemmas \ref{lemma:lgsols} and \ref{lemma:lg2}. We obtain a fundamental 
matrix
\[ \Phi = \left( \begin{array}{llll} U_L & U_R \\
                                     V_L & V_R \end{array}\right), \]
where
\[  U_L = (u_{\phi_1} \ u_{\phi_2})                 
                  = \left( \begin{array}{ll} \phi_1 & \phi_2 \\
                                            \phi'_1 & \phi'_2\end{array}\right), 
\;\;
   V_L =  (v_{\phi_1} \ v_{\phi_2}) 
      = \left( \begin{array}{cc} s\phi'_1-(p\phi''_1)' & s\phi'_2-(p\phi''_2)' 
\\
                                        p\phi''_1 & p\phi''_2 
\end{array}\right), \] 
\[  U_R = (u_{\psi_1} \ u_{\psi_2})
                  = \left( \begin{array}{ll} \psi_1 & \psi_2 \\
                                            \psi'_1 & \psi'_2\end{array}\right), 
\;\;
   V_R = (v_{\psi_1} \ v_{\psi_2}) 
      = \left( \begin{array}{cc} s\psi'_1-(p\psi''_1)' & s\psi'_2-(p\psi''_2)' 
\\
                                        p\psi''_1 & p\psi''_2 
\end{array}\right). \]
In the following, $A^T$ denotes the transpose of the matrix $A$.
\begin{Lemma}\label{lemma:inverse}
Suppose that the boundary conditions at $x=b$ are real.  Then:
\begin{description}
\item[(1)]   $\ds U_L^TV_L - V_L^TU_L = 0 = U_R^TV_R - V_R^TU_R$. 
\item[(2)]  The functions $\psi_1, \psi_2$ can be chosen so that \ \
            $\ds U_R^TV_L - V_R^TU_L = I$ \ (the identity matrix);
\item[(3)]   If $\psi_1, \psi_2$ are chosen to satisfy (2), then
$\ds \Phi^{-1} = \left(\begin{array}{rr} -V_R^T & U_R^T \\
                                      V_L^T & -U_L^T \end{array}\right)$.
\end{description}
\end{Lemma}
{\bf Proof}: (1)  By Lemma \ref{mmlem10}, if $\phi$ and $\psi$
are solutions of $\ell y = \lam y$, then $[\phi,\ov{\psi}](x)$ is constant.

\noindent Since $[\phi_1,\ov{\phi}_2](a) = 0 = [\psi_1,\ov{\psi}_2](b)$, it 
follows
that $[\phi_1,\ov{\phi}_2](x) = 0 = [\psi_1,\ov{\psi}_2](x)$ for all $x \in 
(a,b)$.
Furthermore, $[\phi,\ov{\phi}](x) = 0$ for all functions $\phi(x)$.  This 
implies (1). 
\vskip 5pt

\noindent (2)  Let $S$ denote the 4-dimensional space of solutions of the 
equation
$\ell y = \lam y$.  For \mbox{$\zeta, \eta \in S$}, \ $[\zeta, \ov{\eta}](x) = 
{\rm constant}$.
We shall denote this constant by $[\zeta, \ov{\eta}]$.  The map 
$(\zeta, \eta) \mapsto [\zeta, \ov{\eta}]$ defines a  skew-symmetric,
bilinear form on $S\times S$.  It was shown in the proof of Lemma 
\ref{lemma:lg1} 
that this form is nondegenerate.
\vskip 5pt  

\noindent Let $S_L = {\rm Span}(\phi_1,\phi_2)$, and
$S_R = {\rm Span}(\psi_1,\psi_2)$.  If $ \zeta \in S_L\cap S_R$, then $\zeta$ is 
an
eigenfunction with eigenvalue $\lam$ for the self-adjoint operator $L$.  Since
$\lam \in \Cc \backslash \Rr$, this implies that $\zeta = 0$, and 
\mbox{$S_L\cap S_R = \{0\}$}.
\vskip 10pt

\noindent If $\zeta, \eta \in S_L$, then $[\zeta, \ov{\eta}] = 0$.  Furthermore,
$S_L$ is maximal with respect to this property:  if $\zeta \in S$ and 
$[\zeta, \ov{\eta}] = 0$ for all $\eta \in S_L$, then $\zeta \in S_L$.  (This 
follows
from the fact that if $V$ is a maximal isotropic subspace of a symplectic space 
$W$, 
where $\dim{W} = 2n$, then $\dim{V} = n$.)  $S_R$ is  also maximal with respect
to this property.  This implies that $[\zeta, \ov{\eta}]$
defines a nondegenerate bilinear form on $S_R \times S_L$.  Therefore, 
corresponding
to the basis $\{\phi_1, \phi_2\}$ for $S_L$, there is a dual basis 
$\{\psi_1^o, \psi_2^o\}$ for $S_R$:  $[\psi_i^o, \ov{\phi}_j] = \del_{ij}$ \ 
$(i, j = 1, 2)$.  If we replace $\{\psi_1, \psi_2\}$ by $\{\psi_1^o, 
\psi_2^o\}$,
then $U_R^TV_L - V_R^TU_L = I$.
\vskip 10pt

\noindent (3) From (1) and (2) it follows that 
\[ \left(\begin{array}{rr} -V_R^T & U_R^T \\
                                      V_L^T & -U_L^T \end{array}\right)
   \left(\begin{array}{rr} U_L & U_R \\
                                      V_L & V_R \end{array}\right) =
    \left(\begin{array}{rr} I_2 & 0_2 \\
                                      0_2 & I_2 \end{array}\right) = I_4,  \]
where $I_n$ denotes the $n\times n$ identity matrix, and $0_n$ denotes 
the $n\times n$ zero matrix.  \hfill $\Box$
\vskip 10pt

We will now solve equation (\ref{eq:lg41}) by variation of parameters.  Let
\be z(x) = \Phi(x) \left(\begin{array}{cc} \zeta(x) \\
                                 \eta(x) \end{array}\right), \label{eq:lg42} \ee
where
\[ \zeta(x) = \left(\begin{array}{cc} \zeta_1(x) \\
                                 \zeta_2(x) \end{array}\right), \hm
   \eta(x) = \left(\begin{array}{cc} \eta_1(x) \\
                                 \eta_2(x) \end{array}\right)  \]
are unknown functions.  In the notation preceding equation
(\ref{eq:lg41}), $z = z_y$, where 
\be y = \zeta_1\phi_1 + \zeta_2\phi_2 + \eta_1\psi_1 + \eta_2\psi_2. 
\label{eq:lg43} \ee
Since $y(x)$ must satisfy the boundary conditions at the endpoints,
we require
\be \zeta(b) = 0 = \eta(a).  \label{eq:lg44} \ee
Substituting (\ref{eq:lg42}) into (\ref{eq:lg41}), and using the fact
that the columns of $\Phi$ satisfy the corresponding homogeneous equation, we
obtain
\[ J\Phi \left(\begin{array}{cc} \zeta'(x) \\
                                 \eta'(x) \end{array}\right) = \hat{f}w, \]
or
\[  \left(\begin{array}{cc} \zeta'(x) \\
                                 \eta'(x) \end{array}\right) = 
            \Phi^{-1}\left( \begin{array}{rrrr} 0 \\ 
                                                0 \\                             
                                                        -f \\
                                                0 \end{array} \right)w. \]
By Lemma \ref{lemma:inverse}, if the boundary conditions are real then
\[ \left(\begin{array}{cc} \zeta'(x) \\
                                 \eta'(x) \end{array}\right) =
    \left(\begin{array}{rr} -V_R^T & U_R^T \\
                                      V_L^T & -U_L^T \end{array}\right)
             \left( \begin{array}{rrrr} 0 \\ 
                                        0 \\                                     
                                                -f \\
                                        0 \end{array} \right) =
             \left( \begin{array}{rrrr} -\psi_1 \\ 
                                        -\psi_2 \\                               
                                                  \phi_1 \\
                                         \phi_2 \end{array} \right)fw.  \]
Therefore, using (\ref{eq:lg44}), we have
\[  \zeta(x) = \int_x^b \left(\begin{array}{cc} \psi_1(t) \\
                            \psi_2(t) \end{array}\right) f(t)w(t)\,dt,  \hm                        
             \eta(x) = \int_a^x \left(\begin{array}{cc} \phi_1(t) \\
                            \phi_2(t) \end{array}\right) f(t)w(t)\,dt.  \] 
>From (\ref{eq:lg43}) we see that the solution $y(x)$  of equation 
(\ref{eq:lg41}) is
\[ y(x) = \int_a^bG(x,t)f(t)w(t)\,dt,  \]
where the Green's function $G(x,t)$ \un{\bf for real boundary
conditions} at $x=b$, is 
\be G(x,t) = \left\{\begin{array}{ll} 
    \phi_1(x)\psi_1(t) + \phi_2(x)\psi_2(t) & \mbox{for $a<x<t<b$} \\
    \phi_1(t)\psi_1(x) + \phi_2(t)\psi_2(x) & \mbox{for $a<t<x<b$.}
         \end{array}\right.   \label{eq:green}            \ee
Letting $S = (L-\lam)^{-1}$, we see that for $f \in L^2(a,b;w)$,

\[ (Sf)(x) = (L-\lam)^{-1}f(x) = \int_a^bG(x,t)f(t)w(t)\,dt.  \]
The calculation for complex boundary conditions at $x=b$ is similar.
Lemma \ref{lemma:inverse} remains the same, except that 
\[ U_R^TV_R - V_R^TU_R = K = \left(\begin{array}{rr} 0 & 1 \\
                                                    -1 & 0 \end{array}\right) \]
and
\[ \Phi^{-1} = \left(\begin{array}{cc} -V_R^T-KV_L^T & U_R^T+KU_L^T \\
                                      V_L^T & -U_L^T \end{array}\right).  \]
This implies that
\[  \zeta(x) = \int_x^b \left(\begin{array}{cc} \psi_1(t)+\phi_2(t) \\
                            \psi_2(t)-\phi_1(t) \end{array}\right) f(t)w(t)\,dt, 
 \hm                        
             \eta(x) = \int_a^x \left(\begin{array}{cc} \phi_1(t) \\
                            \phi_2(t) \end{array}\right) f(t)w(t)\,dt,  \]
and the Green's function \un{\bf for complex boundary conditions} at
$x=b$ is
\be G(x,t) = \left\{\begin{array}{ll} 
    \phi_1(x)\psi_1(t) + \phi_2(x)\psi_2(t)
    + \phi_1(x)\phi_2(t) - \phi_2(x)\phi_1(t) & \mbox{for $a<x<t<b$} \\
    \phi_1(t)\psi_1(x) + \phi_2(t)\psi_2(x) & \mbox{for $a<t<x<b$.}
         \end{array}\right.   \label{eq:green1}            \ee
We shall now truncate the interval near $x=b$.  Suppose that the boundary 
conditions
for $L$ at $x=b$ are $[y,\th_1](b) = 0 = [y,\th_2](b)$, where $\th_1$ and 
$\th_2$
are solutions of $\ell y = \lam_0 y$ for some real $\lam_0$.  As before, we 
shall
first carry out the calculation for real boundary conditions, and then indicate 
the
result for complex boundary conditions.

Let $b_j \nearrow b$, and consider the truncated operator $L_j$  defined on 
$(a,b_j]$
with the same boundary conditions as $L$ at $x=a$, and the boundary conditions
$[y,\th_1](b_j) = 0 = [y,\th_2](b_j)$ at $x=b_j$.  The Green's function 
$G_j(x,t)$
will resemble $G(x,t)$ in (\ref{eq:green}), except that $\psi_1$ and $\psi_2$ 
are
replaced by solutions $\psi_1^{(j)}$ and $\psi_2^{(j)}$ which satisfy the 
boundary
conditions for $L_j$ at $x=b_j$.

\begin{Lemma}\label{lemma:lg3} Suppose that $\lambda \in \Cc \backslash \Rr$. 
Then
there are two linearly independent solutions of $\ell y = \lam y$ which satisfy 
the boundary conditions of $L_j$ at $x=b_j$ and are of the form
\be \psi_1^{(j)} = \psi_1 + c_1^{(j)}\phi_1 + c_2^{(j)}\phi_2, \hqs 
     \psi_2^{(j)} = \psi_2 + d_1^{(j)}\phi_1 + d_2^{(j)}\phi_2, \label{eq:102} 
\ee
where
\[ \lim_{j\to \infty} c_1^{(j)} = \lim_{j\to \infty} c_2^{(j)} =
   \lim_{j\to \infty} d_1^{(j)} = \lim_{j\to \infty} d_2^{(j)} = 0. \]
\end{Lemma}

\noindent {\bf Proof}:\,We shall solve for the constants $c_k^{(j)},\ d_k^{(j)}$
in (\ref{eq:102}).  Setting $[\psi_1^{(j)},\th_1](b_j) = 0 = 
[\psi_1^{(j)},\th_2](b_j)$
and $[\psi_2^{(j)},\th_1](b_j) = 0 = [\psi_2^{(j)},\th_2](b_j)$, we obtain the
following equations:

\begin{eqnarray}
c_1^{(j)}\,[\phi_1,\th_1](b_j) + c_2^{(j)}\,[\phi_2,\th_1](b_j) & = & 
-[\psi_1,\th_1](b_j), 
\nonumber  \\
c_1^{(j)}\,[\phi_1,\th_2](b_j) + c_2^{(j)}\,[\phi_2,\th_2](b_j) & = & 
-[\psi_1,\th_2](b_j), 
\nonumber  \\
d_1^{(j)}\,[\phi_1,\th_1](b_j) + d_2^{(j)}\,[\phi_2,\th_1](b_j) & = & 
-[\psi_2,\th_1](b_j), 
\nonumber   \\
d_1^{(j)}\,[\phi_1,\th_2](b_j) + d_2^{(j)}\,[\phi_2,\th_2](b_j) & = & 
-[\psi_2,\th_2](b_j).
              \label{eq:103}  \end{eqnarray}

\noindent  There is no nontrivial solution to the equations
\begin{eqnarray}
c_1\,[\phi_1,\th_1](b_j) + c_2\,[\phi_2,\th_1](b_j) & = & 0, \nonumber \\
c_1\,[\phi_1,\th_2](b_j) + c_2\,[\phi_2,\th_2](b_j) & = & 0, \nonumber
\end{eqnarray}
\noindent since this would imply that $c_1\phi_1 + c_2\phi_2$ is an 
eigenfunction
of the truncated problem with eigenvalue $\lam \in \Cc \sm \Rr$.
Therefore
\[ \Delta(b_j) := \left|\begin{array}{cc} 
                  \left[ \phi_1,\th_1\right](b_j) & 
\left[\phi_2,\th_1\right](b_j) \\
                  \left[ \phi_1,\th_2\right](b_j) & 
\left[\phi_2,\th_2\right](b_j)
                                         \end{array} \right|  
                                         \neq 0, \]
and similarly
\[ \Delta(b) := \left|\begin{array}{cc} \left[\phi_1,\th_1\right](b) & 
\left[\phi_2,\th_1\right](b) \\
                                        \left[\phi_1,\th_2\right](b) & 
\left[\phi_2,\th_2\right](b)
                                          \end{array}\right| 
                                          \neq 0. \]
 Therefore the equations (\ref{eq:103}) have a unique solution:
\[ c_1^{(j)} = \frac{1}{\Delta(b_j)}\left|\begin{array}{cc}
                         \left[\phi_2,\th_1\right](b_j) & 
\left[\psi_1,\th_1\right](b_j)  \\
                         \left[\phi_2,\th_2\right](b_j) & 
\left[\psi_1,\th_2\right](b_j)
                         \end{array}\right|,   \]                         
\[ c_2^{(j)} = \frac{1}{\Delta(b_j)}\left|\begin{array}{cc}
                         \left[\psi_1,\th_1\right](b_j) & 
\left[\phi_1,\th_1\right](b_j)  \\
                         \left[\psi_1,\th_2\right](b_j) & 
\left[\phi_1,\th_2\right](b_j)
                         \end{array}\right|,   \]
and similar formulas for $d_1^{(j)}$ and $d_2^{(j)}$.
Since $\lim_{j \to \infty}[\phi_i,\th_k](b_j) = [\phi_i,\th_k](b)$, and 
$\lim_{j \to \infty}[\psi_i,\th_k](b_j) = [\psi_i,\th_k](b) = 0$,
it follows that
$\ds \lim_{j\to \infty} c_1^{(j)} = \lim_{j\to \infty} c_2^{(j)} =
   \lim_{j\to \infty} d_1^{(j)} = \lim_{j\to \infty} d_2^{(j)} = 0$.  
\hfill $\Box$ 
\vskip 10pt

We now can find the Green's function for the truncated operator $L_j$.
It follows the same pattern as (\ref{eq:green}); the formula is
{\scriptsize
\[ G_j(x,t) =  G(x,t) + c_1^{(j)} \phi_1(x)\phi_1(t) + 
c_2^{(j)}\phi_1(x)\phi_2(t)
+ d_1^{(j)} \phi_2(x)\phi_1(t) + d_2^{(j)}\phi_2(x)\phi_2(t)  
\hspace{5mm} \mbox{for $a<x<t<b_j$.}  \]
\[ G_j(x,t) = G(x,t) + c_1^{(j)} \phi_1(t)\phi_1(x) + 
c_2^{(j)}\phi_1(t)\phi_2(x) 
+ d_1^{(j)} \phi_2(t)\phi_1(x) + d_2^{(j)}\phi_2(t)\phi_2(x)  
\hspace{5mm} \mbox{for $a<t<x<b_j$.}  \] }
The above formula was calculated for the case of real boundary conditions
at $x=b$.  An analogous calculation shows that the formula for complex boundary
conditions is the same, with constants that approach $0$ as $j \rightarrow 
\infty$.

\noindent  As described in section \ref{subsubsection:sisect}, we consider 
$L^2(a,b_j;w) \subset L^2(a,b;w)$, and the splitting
\[ L^2(a,b;w) = L^2(a,b_j;w) \oplus L^2(b_j,b;w). \]
The operator $L_j$ with domain $D(L_j) \subset L^2(a,b_j;w)$ is extended
to $L'_j$ with domain $D(L_j) \oplus L^2(b_j,b;w)$ by the formula
\[ L'_j = L_j \oplus \Th_j,\] 
where $\Th_j$ is the zero operator on $L^2(b_j,b;w)$.  If $P_j$
denotes the projection of $L^2(a,b;w)$ onto $L^2(a,b_j;w)$, then $L'_j = 
L_jP_j$.
Setting $S_j = (L_j-\lam)^{-1}$ and $S'_j = (L'_j-\lam)^{-1}$, we have
\[  (L'_j-\lam)f = \left\{\begin{array}{cc} (L_j-\lam)f & \mbox{for $f\in 
L^2(a,b_j;w)$} \\
                                               -\lam f   & \mbox{for $f\in 
L^2(b_j,b;w)$},
                                               \end{array} \right.  \]
\[ S'_j(f) = (L'_j-\lam)^{-1}f = \left\{\begin{array}{cc} 
                      (L_j-\lam)^{-1}f & \mbox{for $f\in L^2(a,b_j;w)$} \\
                       -\frac{1}{\lam}f & \mbox{for $f\in L^2(b_j,b;w)$},
                                               \end{array} \right.  \]
\[ P_jS'_j = S'_jP_j = S_jP_j.  \]
Let $A = (a,b)\times (a,b), \ \ A_j = (a,b_j)\times (a,b_j)$ and set
\[ G'_j(x,t) = \left\{\begin{array}{ll} 
    G_j(x,t) & \mbox{for $(x,t) \in A_j$} \\
          0  & \mbox{for $(x,t) \in A \sm A_j$.}
         \end{array}\right.               \]
Then
\[ (S_jP_jf)(x) = \int_a^b G'_j(x,t)f(t)w(t)\,dt.   \]
If $L$ has real boundary conditions at $x=b$, then
{\scriptsize
\[ G(x,t) -G'_j(x,t) = \left\{\begin{array}{ll}
 c_1^{(j)} \phi_1(x)\phi_1(t) + c_2^{(j)}\phi_1(x)\phi_2(t) 
 + d_1^{(j)} \phi_2(x)\phi_1(t) + d_2^{(j)}\phi_2(x)\phi_2(t) & \mbox{for 
$a<x<t<b_j$} \\
   c_1^{(j)} \phi_1(t)\phi_1(x) + c_2^{(j)}\phi_1(t)\phi_2(x) 
 + d_1^{(j)} \phi_2(t)\phi_1(x) + d_2^{(j)}\phi_2(t)\phi_2(x) & \mbox{for 
$a<t<x<b_j$} \\
       \phi_1(x)\psi_1(t) + \phi_2(x)\psi_2(t) & \mbox{for $x<t$, \ and \ 
$b_j<t<b$} \\
       \psi_1(x)\phi_1(t) + \psi_2(x)\phi_2(t) & \mbox{for $t<x$, \ and \ 
$b_j<x<b$}.
       \end{array}\right.        \] }
The formula shows that  
$G - G'_j \in L^2(A;w(x)w(t))$ and $G - G'_j \rightarrow 0$ in the
$L^2$ norm as $j \rightarrow \infty$. An analogous calculation shows that the
same is true for complex boundary conditions at $x=b$.
This proves the following theorem.

\begin{Theorem}\label{theorem:halftrunc}
\hspace{1.0in}
\begin{description}
\item[(1)] $S - S_jP_j$ is a Hilbert-Schmidt integral operator. 
\item[(2)] $S - S_jP_j\rightarrow 0$ in the Hilbert-Schmidt norm.
\end{description}
\end{Theorem}

\noindent The theorem implies the following corollary.

\begin{Corollary}\label{corollary:halftrunc}
\hspace{1.0in}
\begin{description}
\item[(1)]  The sequence of truncated operators $L_j$ is spectrally exact
for $L$.
\item[(2)]  $L_j$ and $L$ have the same essential spectrum.
\end{description}
\end{Corollary}

\subsection{The lim-4, lim-4 case}

We now suppose that both endpoints are lim-4, and we shall consider double
truncations.  In this case, $L$ has discrete spectrum.  By translating the
operator if necessary, we may suppose that $0$ is not an eigenvalue.
The eigenvalues of $L$ may be indexed by positive and negative indices:
\[ \cdots \le \lam_{-3} \le \lam_{-2} \le \lam_{-1} < 0 < \lam_1 \le \lam_2 \le
          \lam_3 \le \cdots \le \lam_n \le \cdots   \]
\noindent  It is possible that $L$ is not bounded below, and there are
infinitely many negative eigenvalues.

Let $G(x,t)$ be the Green's function for $L$ at $\lam=0$.  Thus
\[ (L^{-1}f)(x) = \int_a^b G(x,t)f(t)w(t)\,dt.      \] 
There exist functions
$\phi_1,\ \phi_2,\ \psi_1,\ \psi_2 \in L^2(a,b;w)$ which are solutions of
$\ell y = 0$, and such that $\phi_1,\ \phi_2$ provide the boundary conditions
for $L$ at $x=a$, and $\psi_1,\ \psi_2$ provide the boundary conditions at 
$x=b$.
$G(x,t)$ can be expressed in terms of these functions by formulas similar to
those in the previous section.  There are now four cases, depending on whether
the boundary conditions are real or complex at $x=a$ and $x=b$.  We shall not
require the precise formulas.  

Let $a_j \searrow a, \ \ b_j \nearrow b$ and consider the truncated operators
$L_j$ on the intervals $[a_j,b_j]$ with boundary conditions
$[y,\phi_1](a_j) = 0 = [y,\phi_2](a_j)$ and $[y,\psi_1](b_j) = 0 = 
[y,\psi_2](b_j)$.
The solutions $\phi_1^{(j)},\ \phi_2^{(j)}$ of $\ell y = 0$ satisfying the 
boundary
conditions at $x=a_j$ are simply the restrictions $\phi_i^{(j)} = 
\phi_i|[a_j,b_j]$;
and similarly the solutions satisfying the boundary conditions at $x=b_j$ are
$\psi_i^{(j)} = \psi_i|[a_j,b_j]$

Let $S_L = {\rm span}\,(\phi_1,\phi_2)$ and $S_R = {\rm span}\,(\psi_1,\psi_2)$.
$S_L \cap S_R = \{0\}$, since $0$ is not an eigenvalue of $L$.  This implies
that $0$ is not an eigenvalue of $L_j$.  Let $A = (a,b)\times (a,b)$ and
$A_j = (a_j,b_j)\times (a_j,b_j)$.  The Green's function $G_j(x,t)$ for $L_j$ 
at $\lam = 0$
is defined in $A_j$ and coincides with $G(x,t)$ there.

We have the splitting 
\[ L^2(a,b;w) = L^2(a_j,b_j;w) \oplus L^2(a_j,b_j;w)^{\perp}, \] 
\noindent where
\[ L^2(a_j,b_j;w)^{\perp} = L^2(a,a_j;w) \oplus L^2(b_j,b;w). \]
\noindent  Let $S_j = L_j ^{-1}$ and 
\[ S'_j = S_j \oplus \Th_j, \]
\noindent where $\Th_j$ is the zero operator on $L^2(a_j,b_j;w)^{\perp}$.
Let $P_j$ be the projection of $L^2(a,b;w)$ onto $L^2(a_j,b_j;w)$.
Then $P_jS'_j = S'_jP_j = S_jP_j$.  $S_jP_j$ is an integral operator
with kernel
\[ G'_j(x,t) = \left\{\begin{array}{cc} G(x,t) & \mbox{for $(x,t)\in A_j$}, \\
                                            0  & \mbox{for $(x,t)\in A\sm A_j$}.
                                            \end{array}\right.  \]
\noindent It is now clear that
$G - G'_j \in L^2(A;w(x)w(t))$ and $G'_j \rightarrow G$ in the $L^2$
norm.  We have proved the following.

\begin{Theorem}\label{Theorem:trunc}
$S_jP_j \rightarrow S$ in the Hilbert-Schmidt norm.
\end{Theorem}

\noindent  Let the eigenvalues of $L_j$ be:
\[ \cdots \le \lam_{-3}^{(j)} \le \lam_{-2}^{(j)} \le \lam_{-1}^{(j)} < 0 < 
\lam_1^{(j)} 
  \le \lam_2^{(j)} \le  \lam_3^{(j)} \le \cdots \le \lam_n^{(j)} \le \cdots  \]
\noindent  Since $L_j$ is a regular operator, it has only finitely many negative
eigenvalues.  Nevertheless, for any (positive) integer $k$, if $\lam_{-k}$ 
exists,
then for sufficiently large $j$,  $\lam_{-k}^{(j)}$ exists, and 
$\lam_{-k}^{(j)} \rightarrow \lam_{-k}$ \ as \ $j \rightarrow \infty$.

\begin{Theorem}
For $k > 0$, $\lam_k^{(j)} \rightarrow \lam_k$ \ as \ $j \rightarrow \infty$. 
Furthermore, if $\lam_{-k}$ exists, then for sufficiently large $j$, $L_j$ has 
an eigenvalue $\lam_{-k}^{(j)}$, and 
$\lam_{-k}^{(j)} \rightarrow \lam_{-k}$ \ as \ $j \rightarrow \infty$.
\end{Theorem}
\vskip 5pt

\noindent {\bf Proof}: $S$ and $S_jP_j$ are Hilbert-Schmidt operators, 
and $S_jP_j \rightarrow S$ in the Hilbert-Schmidt norm.
These are self-adjoint, compact operators with eigenvalues 
$ \mu_k = 1/{\lam_k}$ \ and \ $ \mu_k^{(j)} = 1/{\lam_k^{(j)}}$,
respectively.  ($S_jP_j$ also has $0$ as an eigenvalue of infinite multiplicity,
but this is not related to $L_j$ or $L$, and can be ignored.)  
The negative eigenvalues are given by a min-max variational
principle, and the positive eigenvalues by a max-min principle.
For $k > 0$:
\[  \mu_k = 
      \max_{M_k}\min_{\stackrel{x\in M_k}{\scriptscriptstyle \|x\|=1}}\langle 
Sx,x\rangle, \]
\[  \mu_{-k} = 
      \min_{M_k}\max_{\stackrel{x\in M_k}{\scriptscriptstyle \|x\|=1}}\langle 
Sx,x\rangle, \]
\noindent where $M_k$ runs through the $k$-dimensional subspaces of
$L^2(a,b;w)$.  The $\mu_k^{(j)}$ are given similarly in terms of
$\langle S_jP_jx,x\rangle$.  The result now follows from the fact that
\[ |\langle(S-S_jP_j)x,x\rangle| \ \ \le \ \ \|S-S_jP_j\| \rightarrow 0, \ \ 
                   {\rm as} \ \ j \rightarrow \infty,   \]
\noindent for $\|x\| = 1$. \hfill $\Box$
\vskip 10pt

\section{One singular endpoint: the lim-2 case}

\setcounter{equation}{0}
\setcounter{Theorem}{0}
\setcounter{Corollary}{0}
\setcounter{Lemma}{0}
\setcounter{Definition}{0}

In this section we consider a problem with one regular endpoint and
one lim-2 singular endpoint. At the regular endpoint, say $x=a$, we
impose two regular self-adjoint boundary conditions. The pre-minimal
domain is then precisely the set $C_{min}$ which we described in
Section \ref{section:Friedrichs} for the case of one singular
endpoint, and the minimal domain is its graph-norm closure. Because
the endpoint $x=b$ is lim-2 the minimal operator possesses only one
self-adjoint extension $L$. If the minimal operator is bounded below,
then $L$ is the Friedrichs extension, and the results of Section
\ref{section:Friedrichs} show how we can obtain spectral exactness
below the essential spectrum; otherwise there is spectrum extending to  
$-\infty$, possibly with gaps. If there are no gaps then spectral
inclusion evidently implies spectral exactness; in the event of gaps, 
we will show that spectral exactness may not be obtained.

The regular approximations $L_j$ to $L$ can still be formed as in Section
\ref{section:Friedrichs}. Any element of the pre-minimal domain -- which
is a core for $L$ in the lim-2 case -- will
satisfy the Dirichlet conditions $y(b_j) = y'(b_j) = 0$ for all
sufficiently large $j$, and so its restriction to the intervals
$[a,b_j]$ will be in the domain of $L_j$ for all sufficiently large $j$.
This means that the spectra of the $L_j$ give spectrally inclusive
approximations to the spectrum of $L$. In particular, we have
the following result.
\begin{Proposition}\label{proposition:unbounded}
Suppose that the minimal operator is unbounded below and, for each
fixed integer $k\geq 0$, let $\mu_k^{(j)}$ be the $k$th eigenvalue
of $L_j$. Then
\[ \lim_{j\rightarrow\infty}\mu_k^{(j)} = -\infty. \]
\end{Proposition}
It is also not difficult to see that if there are gaps in the essential
spectrum of $L$ then the $b_j$ may be chosen to ensure that for some point
$\lambda^*$ in one of these gaps, there is some $k$ such that
\[ \mu_k^{(j)} = \lambda^*. \]
To see how this may be achieved let $\beta$ be less than $b$ and let
$\mu_k(\beta)$ be the $k$th eigenvalue for the problem on $[a,\beta]$
with Dirichlet conditions at $\beta$. By taking $k$ sufficiently large
we can ensure that $\mu_k(\beta) > \lambda^*$; with $k$ now fixed
we can let $\beta$ increase towards $b$ until $\mu_k(\beta)$ attains
the value $\lambda^*$. Thus we obtain the following result.

\begin{Proposition}\label{prop:spurious}
Suppose that $L$ is unbounded below and suppose that $\lambda^*$
is not a spectral point of $L$. Then the approximating operators
$L_j$ may be constructed so that $\lambda^*$ is a spectral point
of {\bf every} $L_j$.
\end{Proposition}

The approximation of essential spectrum evidently requires something
more than eigenvalues: it requires the spectral function or the
Titchmarsh-Weyl $M(\lambda)$ matrix. Computing the spectral function
is a difficult problem even in the second order case. For more 
information about ways of computing the $M(\lambda)$ matrix in the
fourth order case, see \cite{Brownetal}.

\section{One singular endpoint: the lim-3 case}\label{section:onelim3}

\setcounter{equation}{0}
\setcounter{Theorem}{0}
\setcounter{Corollary}{0}
\setcounter{Lemma}{0}
\setcounter{Definition}{0}

We shall now consider a problem where the endpoint $x=a$ is regular, and 
$x=b$ is lim-3 singular.  A self-adjoint extension $L$ of $L_{min}$ will
have boundary conditions of the form 
\be A_1u_y(a) + A_2v_y(a) = 0,  \label{pr2} \ee
\be [y,\psi](b) = 0,  \label{pr3} \ee
where $\psi$ is a real function in $D_{max} \sm D_{min}$.  
In the following, we will need
to compare $L$ with another self-adjoint extension $L_F$.  Let $L_0$ be the
extension of $L_{min}$ with boundary conditions (\ref{pr2}) at $x=a$.
If $L$ is bounded below, then so is $L_0$, and therefore $L_0$ has a
Friedrichs extension $L_F$.  This operator will play an important role
in the following. 
 
Suppose that $L$ is bounded below, with a number of eigenvalues strictly
below the essential spectrum -- say $\lambda_{0},
\ldots,\lambda_{n}$. {\bf Note}: we are not assuming that these are
the only eigenvalues below the essential spectrum. We shall construct
a sequence of regular operators $L_j$  on truncated intervals $[a,b_j)$,
such that for each $k$ between $0$ and $n$, the $k$th eigenvalue of $L_j$
converges to the $k$th eigenvalue of $L$ as $b_j \nearrow b$. This will 
require a careful choice of boundary conditions at the endpoints $b_j$, 
which we now describe.
\vskip 5pt

\noindent We shall denote by
\[ \left(\begin{array}{l} U_L(x) \\ V_L(x) \end{array}\right) \]
a $4\times 2$ fundamental matrix solution of the Hamiltonian form of the
Sturm-Liouville equation which satisfies the initial conditions 
(\ref{eq:lgintro})
at $x=a$. The dependence on $\lambda$ is suppressed in this notation.
Associated with this fundamental matrix is the matrix
\be W_L(x) = V_L(x)U_L^{-1}(x), \label{eq:mm13} \ee
which is defined except at a finite number of points in $(a,b)$
(assuming that $\lam$ lies below the essential spectrum). 
Suppose that, in the Hamiltonian
formulation, the boundary conditions for $L_j$ at $x=b_j$ have the form
\be v = W_R u, \label{eq:mm5} \ee
where $W_R$ is a real, symmetric matrix:
\be W_R = \left(\begin{array}{ll} \kappa & \mu \\ \mu & \nu 
 \end{array}\right). \label{eq:mmrev1} \ee
By Theorem \ref{Theorem:intro}, the number of eigenvalues of
$L_j$ which are strictly less than $\lam$ is
\be N(b_j,\lambda) = \delta_{L}(b_j,\lambda) + \sigma(b_j,\lambda), 
\label{eq:mm10} \ee
where
\[ \delta_L(b_j,\lambda) = \sum_{a < x < b_j}\mbox{(Rank deficiency of 
$U_{L}(x)$)}. \]
If $\det U_L(b_j,\lam) \neq 0$, then
\[   \sigma(b_j,\lambda) = \nu_{\#}(W_L(b_j)-W_R),  \]
where $\nu_{\#}(W)$ is the number of negative eigenvalues of $W$.
We shall investigate the possible
values for $\sigma(b_j,\lambda)$ corresponding to different choices of the 
matrix
$W_R$ which defines the boundary condition at $x=b_j$.  In particular, we shall 
show that $\sigma=1$ is always possible. We shall use the notation of 
(\ref{eq:mmqd}) for quasiderivatives. We require the vectors
\[ u_{\psi} = \left(\begin{array}{l} \psi^{[0]} \\ \psi^{[1]} 
\end{array}\right), \;\;\;
   v_{\psi} = \left(\begin{array}{l} \psi^{[3]} \\ \psi^{[2]} 
\end{array}\right). \]
Because the boundary condition defined by (\ref{eq:mm5}) must include
the boundary condition defined by (\ref{pr3}), we must have
\be v_{\psi} = W_R u_{\psi}. \label{eq:mm6} \ee
If we assume that $\psi^{[1]}$ is non-zero, then we can rearrange
this formula to obtain $\mu$ and $\nu$ in terms of the unknown
$\kappa$ and the known $\psi^{[0]}$, $\psi^{[1]}$, $\psi^{[2]}$ and 
$\psi^{[3]}$:
\[ \mu = (\psi^{[3]}-\kappa\psi^{[0]})/\psi^{[1]}, \;\;\;
   \nu = \frac{\psi^{[2]}}{\psi^{[1]}} - 
    \frac{\psi^{[0]}}{(\psi^{[1]})^2}(\psi^{[3]}-\kappa\psi^{[0]}).\]   
To determine the number of negative eigenvalues of $W_L-W_R$ we require
the trace and determinant of this matrix. Expressing these in terms of
$\kappa$ we obtain
\be \mbox{trace}(W_L-W_R) = \mbox{trace}(W_L) +
\frac{(\psi^{[0]}\psi^{[3]}-\psi^{[1]}\psi^{[2]})}{(\psi^{[1]})^2}
 - \kappa\left(1+\frac{(\psi^{[0]})^{2}}{(\psi^{[1]})^{2}}\right). 
\label{eq:mm7} \ee
\be \mbox{det}(W_L-W_R) =  C + 
\frac{\kappa}{(\psi^{[1]})^2}(u_{\psi}^{T}v_{\psi}
    - u_{\psi}^{T}W_L u_{\psi}).  \label{eq:mm8} \ee
The constant $C$ in (\ref{eq:mm8}) is given by 
\be C = \mbox{det}(W_L) + \frac{1}{(\psi^{[1]})^2}\left(k(\psi^{[0]}\psi^{[3]}
 -\psi^{[1]}\psi^{[2]})
 - (\psi^{[3]})^2 + (m+\ov{m})\psi^{[1]}\psi^{[3]}\right),  \label{eq:mm9} \ee
where $k$ and $m$ are, respectively, the (1,1) and (1,2) terms
of the Hermitian matrix $W_L$.

Similarly, if we assume that $\psi^{[0]}$ is non-zero, then we can rearrange
this formula to obtain $\mu$ and $\kappa$ in terms of the unknown
$\nu$ and the known $\psi^{[0]}$, $\psi^{[1]}$, $\psi^{[2]}$ and $\psi^{[3]}$.  
We then
obtain the following formulas for the trace and determinant of $W_L-W_R$:
\be \mbox{trace}(W_L-W_R) = \mbox{trace}(W_L) +
\frac{(\psi^{[1]}\psi^{[2]}-\psi^{[0]}\psi^{[3]})}{(\psi^{[0]})^2}
 - \nu\left(1+\frac{(\psi^{[1]})^{2}}{(\psi^{[0]})^{2}}\right), \label{eq:lg7} 
\ee
\be \mbox{det}(W_L-W_R) =  D + \frac{\nu}{(\psi^{[0]})^2}(u_{\psi}^{T}v_{\psi}
    - u_{\psi}^{T}W_L u_{\psi}),  \label{eq:lg8} \ee
where the constant $D$ in (\ref{eq:lg8}) is given by 
\be D = \mbox{det}(W_L) + \frac{1}{(\psi^{[0]})^2}\left(n(\psi^{[1]}\psi^{[2]}
 -\psi^{[0]}\psi^{[3]})
 - (\psi^{[2]})^2 + (m+\ov{m})\psi^{[0]}\psi^{[2]}\right).  \label{eq:lg9}    
\ee
Here $n$ is the (2,2) term of $W_L$.

We shall now prove a number of results which indicate how the $L_j$ can be 
constructed to obtain spectral exactness. The following lemma shows that for 
a given $\lam_*$ below the
essential spectrum, we can find a nearby value $\lam$ so that 
${\rm det}\,U_L(b_j,\lam) \neq 0$\ \ and\ \ $\sig(b_j,\lam) = 1$.
Recall that $L_F$ has been defined as the Friedrichs extension of the
operator $L_0$, which has boundary conditions (\ref{pr2}).
\vskip 10pt

\begin{Lemma}\label{lemma:sigma}
Suppose that

\begin{description}
\item[(a)]  The operator $L$ is bounded below;
\item[(b)]  The interval $[\lam_*,\lam_*+\eps]$ lies below the essential
 spectrum, and contains no eigenvalue of the Friedrichs extension $L_F$.
\end{description}
Then
\begin{description}
\item[(1)]  There exists $\bet_0 \in [a,b)$ such that for any $\bet \in 
[\bet_0,b)$
            and $\lam \in [\lam_*,\lam_*+\eps]$, ${\rm det}\,U_L(\bet,\lam) \neq 
0$;
\item[(2)]  If $\bet \in [\bet_0,b)$ and $u_{\psi}(\bet) \ne 0$, then there
            exists  $\lam \in [\lam_*,\lam_*+\eps]$ and a 2$\times$2 real, 
symmetric
            matrix  $W_R$ such that:                        
     \begin{description} 
     \item[(i)]  $v_{\psi}(\bet) = W_Ru_{\psi}(\bet)$;
     \item[(ii)]  $\sig(\bet,\lam) = 1$ for the truncated eigenvalue problem
\[\left\{ \begin{array}{cc}
           \ell(y) = \lam y, \ \ a < x < \bet, \\
           A_1u(a) +A_2v(a) = 0, \\
           v(\bet) = W_Ru(\bet).
           \end{array}
\right.  \]
       \end{description}
\end{description}
\end{Lemma}
\vskip 10pt

\noindent {\bf Proof:}\, {\bf (1)}  For $a < \bet < b$, consider the eigenvalue 
 problem
\[\left\{ \begin{array}{cc}
           \ell(y) = \lam y, \ \ a < x < \bet, \\
           A_1u(a) +A_2v(a) = 0, \\
           u(\bet) = 0.
           \end{array}
\right.  \]
We shall denote this eigenvalue problem by EP($\bet$).  Let $\mu_k$ denote
the Friedrichs eigenvalues, and $\mu_k(\bet)$ the eigenvalues of EP($\bet$).  By
Theorem \ref{Theorem:FApprox} $\mu_k(\bet) \searrow \mu_k$ \ as \ 
$\bet \nearrow b$.
\vskip 5pt

\noindent Since $\lam_* + \eps$ lies below the essential spectrum, there are 
only finitely
many Friedrichs eigenvalues (perhaps none) which are less than $\lam_*+\eps$.
\vskip 5pt

\noindent If there are no Friedrichs eigenvalues below $\lam_*+\eps$, then \ 
$\mu_k(\bet) > \lam_*+\eps$ \ for all \ $\bet \in [a,b)$ \ and all indices $k$.
In this case, we take $\bet_0 = a$.
\vskip 5pt

\noindent If there are some Friedrichs eigenvalues below $\lam_*+\eps$, let 
$\mu_N$ 
be the largest one.  Since $[\lam_*,\lam_* + \eps]$ contains no Friedrichs 
eigenvalues,  
$\mu_N < \lam_*$.  This implies that there is $\bet_0 \in [a,b)$ such that 
$\mu_N(\bet)<\lam_*$ for all $\bet \in [\bet_0,b)$.  Furthermore, if $\mu_{N+1}$ 
exists (below the essential spectrum) then 
$\mu_{N+1}(\bet)>\mu_{N+1}>\lam_*+\eps$ \ 
for all $\bet \in (a,b)$.  If $\mu_{N+1}$ does not exist below the essential 
spectrum, 
then $\mu_{N+1}(\bet)$ does not lie below the essential spectrum for any $\bet 
\in (a,b)$.
Thus no eigenvalue of EP($\bet$) lies in $[\lam_*,\lam_*+\eps]$ for $\bet \in 
[\bet_0,b)$.
But $\lam$ is an eigenvalue of EP($\bet$) if and only if ${\rm 
det}\,U_L(\bet,\lam) = 0$.
This shows that ${\rm det}\,U_L(\bet,\lam) \neq 0$ for all $\lam \in 
[\lam_*,\lam_*+\eps]$
and $\bet \in [\bet_0,b)$.
\vskip 10pt

\noindent {\bf (2)}  For $\bet \in [\bet_0,b)$, and $\lam \in 
[\lam_*,\lam_*+\eps]$, 
${\rm det}\,U_L(\bet,\lam) \neq 0$.  Therefore $W_L(\bet,\lam) = 
V_L(\bet,\lam)U_L(\bet,\lam)^{-1}$ is defined.  For the eigenvalue problem in
{\bf (2)(ii)}, $\sig(\bet,\lam)$ is the number of negative eigenvalues of 
$W_L(\bet,\lam) - W_R$. Now $\sig(\bet,\lam) = 1$ if 

\be {\rm det}\left(W_L(\bet,\lam) - W_R\right) < 0. \label{eq:lg10} \ee  

\noindent The matrix $W_R$ is constructed as indicated in the calculations 
preceding
this lemma. (The construction guarantees that $v_{\psi}(\bet) = 
W_Ru_{\psi}(\bet)$.)  
Since $u_{\psi}(\bet) \ne 0$, ${\rm det}\left(W_L(\bet,\lam) - W_R\right)$ 
can be calculated by one of the formulas (\ref{eq:mm8}) if $\psi^{[1]}(\bet) \ne 
0$, 
or (\ref{eq:lg8}) if $\psi^{[0]}(\bet) \ne 0$.  We can then force the inequality 
(\ref{eq:lg10}) by an appropriate choice of $\kap$ or $\nu$, provided that
\be u_{\psi}^Tv_{\psi} - u_{\psi}^TW_L(\bet,\lam)u_{\psi} \ne 0, \label{eq:lg11} 
\ee
where $u_{\psi}$ and $v_{\psi}$ are evaluated at $\bet$.
\vskip 5pt

\noindent  We claim that the inequality (\ref{eq:lg11}) is satisfied for some
$\lam \in [\lam_*,\lam_*+\eps]$.  This follows from the fact (shown in 
Greenberg \cite{kn:LG1}) that $W_L(\bet,\lam)$ is a strictly decreasing
matrix function of $\lam$ in any $\lam$-interval containing no zeros of
${\rm det}\,U_L(\bet,\lam)$.  Thus if 
$u_{\psi}^Tv_{\psi} - u_{\psi}^TW_L(\bet,\lam_*)u_{\psi} = 0$, and if $\lam_*$
is increased slighty to $\lam_{**}$, then the inequality (\ref{eq:lg11})
will be satisfied for $\lam = \lam_{**}$.  (Here we have again used the
fact that $u_{\psi}(\bet) \ne 0$.)  \hfill $\Box$
\vskip 10pt

\noindent \un{\bf Remarks.} \ \  Suppose that $L$ is a self-adjoint extension 
of $L_0$ with domain $D(L) = D(L_0) \oplus {\rm Span}\,(\psi)$.
\vskip 5pt

\noindent (1) The function $\psi$ cannot have compact support.  Thus the
assumption $u_{\psi}(\bet) \ne 0$ (in part {\bf (2)(ii)} of the preceding lemma) 
is 
satisfied on some sequence $b_j \nearrow b$.
\vskip 5pt

\noindent (2)  Suppose that there is a sequence $b_j \nearrow b$ such that
$u_{\psi}(b_j) = 0$, for all $j$.  Then the eigenvalues of $L$ below the
essential spectrum coincide with those of the Friedrichs extension $L_F$.  These 
may
be approximated by regular truncated problems with Dirichlet boundary conditions
at the points $b_j$. 

\begin{Theorem} \label{Theorem:interlacing}
Suppose that the operator $L$ is bounded below, and has eigenvalues
$\lam_0, \ldots,\lam_n$ below the essential spectrum.  Suppose also
that the Friedrichs extension $L_F$ of $L_0$ has eigenvalues
$\mu_0, \ldots, \mu_n$ below the essential spectrum.  {\rm [}{\bf Note}:
we are not assuming that these are the only eigenvalues below the essential 
spectrum.{\rm ]}  Then the following inequalities hold between the eigenvalues 
of $L$ and $L_F$:
\be \lam_0 \leq \mu_0, \ \ \ \mu_{k-1} \leq \lambda_k \leq \mu_{k}, \;\;\; 
 k = 1,\ldots, n.  \label{eq:mm10b}  \ee
\end{Theorem}

\noindent {\bf Proof}:\, The inequality $\lambda_k \leq \mu_k$ is well-known 
(see, e.g., Dunford and Schwartz 
\cite[Problem D2, p.1544]{kn:Dunford}). For the inequality
$\mu_{k-1}\leq \lambda_k$ we proceed
by contradiction. Suppose that $\lambda_k < \mu_{k-1}$, and fix 
$\lambda_*$ and $\eps$ so that 
$[\lam_*,\lam_* +\eps] \subset (\lambda_k,\mu_{k-1})$. By Lemma 
\ref{lemma:sigma}
we can choose a
sequence of points $b_j$ converging to $b$ and a set of associated boundary 
conditions at $b_j$ which give $\sigma(b_j,\nu^{(j)}) = 1$ for some 
$\nu^{(j)} \in [\lam_*,\lam_*+\eps]$. Let $\mu_i(b_j)$ be the eigenvalues of the 
truncated problem on $[a,b_j]$ with Dirichlet boundary conditions at $x=b_j$.
By Lemma \ref{lemma:lf3} we know that for all sufficiently large $j$ we 
have $\lam_* + \eps < \mu_{k-1} < \mu_{k-1}(b_j)$.  Since $\delta_L(b_j,\lam)$
is the number of eigenvalues $\mu_i(b_j)$ less than $\lam$,
this implies that for $\lam \in [\lam_*,\lam_* +\eps], \ \delta_L(b_j,\lam) \leq 
k-1$ 
in (\ref{eq:mm10}).   Since $\sigma(b_j,\nu^{(j)}) = 1$, this means that
$N(b_j,\nu^{(j)}) = \delta_L(b_j,\nu^{(j)}) + \sigma(b_j,\nu^{(j)}) \leq k$.  
Thus the truncated
operator $L_j$ has at most $k$ eigenvalues less than $\nu^{(j)}$, and therefore
at most $k$ eigenvalues less than $\lam_*$.  But the singular operator $L$ has
$k+1$ eigenvalues less than $\lam_*$.  This violates  the spectral inclusion 
guaranteed by Theorem \ref{Theorem:SRCSI}. \hfill $\Box$  
\vspace{2mm}

\noindent If $\lam_n < \mu_n$, then we can prove the following sharper
version of Lemma \ref{lemma:sigma}\,{\bf (2)}.

\begin{Lemma}\label{lemma:sigma1}
Suppose that
\begin{description}
\item[(a)]  $L$ is bounded below;
\item[(b)]  $\lam_n < \lam_* < \mu_n$.
\end{description}
Then there exists $\gam_0 \in [a,b)$ such that if
$\bet \in [\gam_0,b)$ and $u_{\psi}(\bet) \neq 0$, there is
a 2$\times$2 real, symmetric matrix $W_R$ such that
\begin{description} 
\item[(1)]  $v_{\psi}(\bet) = W_Ru_{\psi}(\bet)$;
\item[(2)]  $\sig(\bet,\lam_*) = 1$ for the truncated eigenvalue problem
\[\left\{ \begin{array}{cc}
           \ell(y) = \lam y, \ \ a < x < \bet, \\
           A_1u(a) +A_2v(a) = 0, \\
           v(\bet) = W_Ru(\bet).
           \end{array}
\right.  \]
\end{description}
\end{Lemma}
\vskip 10pt

\noindent {\bf Proof:}\, By Lemma \ref{lemma:sigma}, there is $\bet_0$ such 
that if $\bet \in [\bet_0,b)$ then ${\rm det}\,U_L(\bet,\lam_*) \neq 0$.
Therefore, for the eigenvalue problem in {\bf (2)}, $\sig(\bet,\lam_*)$ 
is the number of negative eigenvalues of \mbox{$W_L(\bet,\lam_*) - W_R$}, so  
$\sig(\bet,\lam_*) = 1$ if 
\be {\rm det}\left(W_L(\bet,\lam_*) - W_R\right) < 0.  \label{eq:lg12} \ee
We construct $W_R$ by the calculations preceding Lemma \ref{lemma:sigma}. 
The trace and determinant of \hfill\linebreak 
\mbox{$W_L(\bet,\lam_*) - W_R$} are given
by the formulas (\ref{eq:mm7}), (\ref{eq:mm8}), (\ref{eq:lg7}), (\ref{eq:lg8}).
The inequality (\ref{eq:lg12}) will be forced by an appropriate choice of
$\kap$ or $\nu$ in (\ref{eq:mmrev1}) if
\be u_{\psi}^Tv_{\psi} - u_{\psi}^TW_L(\bet,\lam)u_{\psi} \ne 0. \label{eq:lg13} 
\ee
If the left hand side of (\ref{eq:lg13}) is zero, then 
${\rm det}\left(W_L(\bet,\lam_*) - W_R\right)$ equals $C = C(\bet)$  in 
(\ref{eq:mm9}) if $\psi^{[1]}(\bet) \neq 0$, or it equals $D = D(\bet)$ in
(\ref{eq:lg9}) if $\psi^{[0]}(\bet) \neq 0$.  If these coefficients are 
negative,
then  $\sig(\bet,\lam_*) = 1$.  To prove the existence of $\gam_0 \in 
[\bet_0,b)$
with the properties stated in the lemma, we shall argue by contradiction.
Suppose that there is a sequence $b_j \nearrow b$ with the following properties:
\begin{description}
\item[(i)]  $u_{\psi}(b_j) \neq 0$;
\item[(ii)] $ u_{\psi}(b_j)^Tv_{\psi}(b_j) -
                    u_{\psi}(b_j)^TW_L(b_j,\lam)u_{\psi}(b_j) = 0$;
\item[(iii)]  $C(b_j) \geq 0$ if $\psi^{[1]}(b_j) \neq 0$; \ \ \
              $D(b_j) \geq 0$ if $\psi^{[0]}(b_j) \neq 0$.
\end{description}
Note that $C(b_j)$ and $D(b_j)$ do not involve $\kap$ or $\nu$. 
Therefore ${\rm det}\left(W_L(b_j,\lam_*) - W_R\right)$ is independent of
$\kap$ or $\nu$, and is nonnegative.  By an appropriate choice of $\kap$ 
or $\nu$, we can construct $W_R$ so that
${\rm trace}\,\left(W_L(b_j,\lam_*) - W_R\right) > 0$.  This implies that
$\sig(b_j,\lam_*) = 0$ for the eigenvalue problem
\[\left\{ \begin{array}{cc}
           \ell(y) = \lam y, \ \ a < x < b_j, \\
           A_1u(a) +A_2v(a) = 0, \\
           v(b_j) = W_Ru(b_j).
           \end{array}
\right.  \]
We shall denote this truncated eigenvalue problem by EP($b_j$).
\vskip 5pt

\noindent  Let $\mu_k(b_j)$ denote the eigenvalues of the truncated problem
on $[a,b_j]$ with Dirichlet boundary conditions at $x = b_j$.
By Lemma \ref{lemma:lf3}, for sufficiently large $j$, 
\[ \lam_n < \lam_* <  \mu_n < \mu_n(b_j).  \]
Since $\delta_L(b_j,\lam_*)$ is the number of eigenvalues
$\mu_k(b_j)$  less than $\lam_*$,
this implies that $\delta_L(b_j,\lam_*) \leq n$, and 
\[  N(b_j,\lam_*) = \delta_L(b_j,\lam_*) + \sig(b_j,\lam_*) \leq n.  \]
Thus the approximating problems EP($b_j$) have at most $n$
eigenvalues less than $\lam_*$, while the singular problem has $n+1$
such eigenvalues. This contradicts spectral inclusion. \hfill $\Box$

\begin{Theorem}\label{Theorem:1}
Suppose that the hypotheses of Theorem \ref{Theorem:interlacing} hold, and 
suppose 
also that $\lambda_{n} < \mu_{n}$. Then we can construct a sequence of
operators $L_j$ on intervals $[a,b_j]$ whose eigenvalues 
$\lambda_{k}(b_j)$ converge, for each $0 \leq k \leq n$, to the corresponding
eigenvalues $\lambda_{k}$ of $L$:
\[ \lim_{j\rightarrow\infty}\lambda_{k}(b_j) = \lambda_k, \;\;\; k = 
0,1,\ldots,n. \]
\end{Theorem}

\noindent {\bf Proof}:\, Fix $\lam_* \in (\lam_n,\mu_n)$, and let $b_j \nearrow 
b$.
By the previous lemma, we can construct truncated eigenvalue problems EP($b_j$):
\[\left\{ \begin{array}{cc}
           \ell(y) = \lam y, \ \ a < x < b_j, \\
           A_1u(a) +A_2v(a) = 0, \\
           v(b_j) = W_ju(b_j),
           \end{array}
\right.  \]
such that $\sig(b_j,\lam_*) = 1$.  Let $\mu_k(b_j)$ denote the eigenvalues
of the truncated problem on $[a,b_j]$ with Dirichlet boundary conditions at 
$x=b_j$. 
For sufficiently large $j$,
\[ \mu_{n-1} < \mu_{n-1}(b_j) < \lam_* < \mu_n < \mu_n(b_j).  \]
This implies that $\delta_L(b_j,\lam_*) = n$, and
\be  N(b_j,\lam_*) = \delta_L(b_j,\lam_*) + \sig(b_j,\lam_*) = n+1.  
\label{eq:N} \ee
Thus the approximating problems EP($b_j$) have $n+1$ eigenvalues less than
$\lam_*$, and the same is true for the singular problem.  Note that 
$\lam_{n+1}(b_j)$
cannot converge to any of the singular eigenvalues less than $\lam_*$, since
equation (\ref{eq:N}) implies that
$\lam_* \leq \lam_{n+1}(b_j)$.  The result now follows from spectral inclusion.
\hfill $\Box$
\vspace{2mm}

This theorem shows that it is possible to construct a sequence of
regular fourth order Sturm-Liouville problems to approximate the eigenvalues
below the essential spectrum
of a lim-3 singular problem, with each such eigenvalue being approximated
by the corresponding regular eigenvalues of the same index. In other words,
we can avoid having a situation in which, say, the eigenvalues $\lambda_0(b_j)
\rightarrow -\infty$
as $j\rightarrow \infty$.  This would be an undesirable phenomenon since
it essentially means that the regular approximating problems possess a 
spurious eigenvalue which is not approximating anything in the spectrum of
the problem which interests us. However, the implementation of the theorem by
a numerical procedure requires that we find a point
$\lambda_*$ in $(\lambda_n,\mu_n)$ for some $n$.   
We may not know enough about the spectrum in advance to be able to choose such 
a point.  The following proposition can be useful for this purpose.

\begin{Proposition} Suppose that the hypotheses of Theorem \ref{Theorem:1}
hold.  In particular, assume that $\lam_n < \mu_n$.
\begin{description}
\item[(i)] Suppose that for some $k\leq n$, $\lambda_k < \mu_k$. Then for
 each $\lambda \in (\lambda_k,\mu_k)$, the minimum value of
 $\sigma(b_j,\lambda)$ is 1 for all sufficiently large $j$.
\item[(ii)] Suppose that for some $k\leq n$, $\mu_{k-1} < \lambda_k$.
 Then for each $\lambda \in (\mu_{k-1},\lambda_k)$ the minimum
 value of $\sigma(b_j,\lambda)$ is 0 for all sufficiently large
 $j$.
\item[(iii)] Suppose that $\lambda < \lambda_0$. Then the minimum
 value of  $\sigma(b_j,\lambda)$ is 0 for all sufficiently large
 $j$.
\end{description}
\end{Proposition}

\noindent {\bf Proof}:\, The proofs of the three parts of this
theorem are all quite similar.  Let $L_j$ indicate the approximating
truncated operators.   From Corollary \ref{Corollary:spec.inclus}
we know that we can choose the $L_j$ to achieve spectrally inclusive
eigenvalue convergence. In Case {\bf (i)}, if $\lambda\in (\lambda_k,\mu_k)$
then $L_j$ has at least $k+1$ eigenvalues less than $\lambda$
for all sufficiently large $j$. The number of eigenvalues of
$L_j$ which are less than $\lambda$ is also given by (\ref{eq:mm10})
in which $\delta_{L}(b_j,\lambda)$ is exactly equal to $k$.  This means that
$\sigma(b_j,\lambda)\geq 1$. However we also showed above that we can choose
the matrix $W_R$ to achieve $\sigma(b_j,\lambda)=1$, so this is the
minimum value. For Case {\bf (ii)}, we exploit Theorem 
\ref{Theorem:1} to assert that we can choose the sequence $L_j$ to be
spectrally exact. In this case we know that for all sufficiently
large $j$, each $L_j$ will have precisely $k$ eigenvalues less
than $\lambda$ because $\lambda_{k-1} \leq \mu_{k-1} < \lambda < \lambda_k$.
Again we know that for all sufficiently large $j$, we have
$\delta_{L}(b_j,\lambda) = k$ for all sufficiently large $j$. Thus
from (\ref{eq:mm10}), this spectrally exact sequence $L_j$ must
be giving us $\sigma(b_j,\lambda)=0$ for all sufficiently large $j$,
which is clearly the minimum value that $\sigma$ can have.
Finally, the proof of Case {\bf (iii)} is virtually identical to
the proof of Case {\bf (ii)}. \hfill $\Box$
\vspace{2mm}

This Lemma gives us a prescription for obtaining eigenvalue approximations.
Suppose we want to compute, for each $\lambda$ below the essential spectrum 
of $L$, the number of eigenvalues of $L$ which are less than $\lambda$. Then all
we need to do is set up a regular approximation over a truncated interval
$[a,b_j]$, with $b_j$ sufficiently close to $b$, and choose our boundary
condition (\ref{eq:mm5}) to minimize $\sigma(b_j,\lambda)$ {\em for this
value of $\lambda$}. $N(b_j,\lambda)$, which we can compute if we can solve
a regular fourth-order Sturm-Liouville problem, will be the eigenvalue
count which we seek.

\section{Two singular endpoints}

\setcounter{equation}{0}
\setcounter{Theorem}{0}
\setcounter{Corollary}{0}
\setcounter{Lemma}{0}
\setcounter{Definition}{0}

\subsection{The lim-2, lim-2 case}
In this case, for a problem posed over $(a,b)$, the preminimal domain
is the set of functions with compact support in $(a,b)$ and the minimal
domain is its closure in the graph norm. However, in this case the
minimal operator is self-adjoint. Thus, if it is bounded below, it 
is its own Friedrichs extension, and the results of Section 
\ref{section:Friedrichs} apply. If it is not bounded below then it is
easy to see, by analogy with Proposition \ref{prop:spurious}, that
if there are gaps in the spectrum then a sequence of approximating 
truncated problems may be constructed in such a way that all will
have an eigenvalue at some fixed point in a gap. Further information
about the spectrum in such cases would therefore require approximations
to the spectral function.

\subsection{The lim-2, lim-3 case}
We consider now a problem posed over an interval $(a,b)$ in which
$x=a$ is a lim-2 endpoint and $x=b$ is lim-3. At $x=a$ we have no
boundary conditions, while at $x=b$ we have a single boundary 
condition of the form $[y,\psi](b)=0$.  We may suppose that $\psi$
is real.  We shall assume that the minimal operator is bounded below.

If the boundary condition function $\psi$ yields the Friedrichs
extension $L_F$ then we have no more work to do: we have already seen how
to obtain spectral exactness in this case. Other cases require
a much more delicate treatment. We shall continue to denote the
eigenvalues of $L_F$ by $\mu_k$; the eigenvalues subject to $[y,\psi](b)=0$ 
will be dnoted by $\lambda_k$.

\begin{Theorem}\label{Theorem:lim2lim3}
Suppose that there is some $n$ such that $\mu_n$ lies strictly
below any essential spectrum and suppose that $\lambda_n 
< \mu_n$. Let $a_j$ be a sequence such that $a_j\searrow a$
as $j\nearrow\infty$, and suppose that Dirichlet boundary conditions
are imposed at $x=a_j$.  Then there exists a sequence $b_j$ with
$b_j\nearrow b$, and a set of boundary conditions at each $b_j$
including the boundary condition $[y,\psi](b_j)=0$, such
that the eigenvalues $\lambda_k^{(j)}$ of the associated
truncated eigenproblems satisfy
\be \lim_{j\rightarrow\infty}\lambda_k^{(j)} = \lambda_k,
 \;\;\; k=0,\ldots,n.
\label{eq:mm24b} 
\ee
\end{Theorem}
\noindent {\bf Proof}:\, Let $\lambda_k(a_j,b)$ denote the
$k$th eigenvalue of the problem with one singular endpoint
at $x=b$, with boundary conditions $y(a_j)=0=y'(a_j)$,
and $[y,\psi](b)=0$. Whatever point $b_j$ we choose and 
whatever boundary condition we impose there, we shall
always have
\be |\lambda_k^{(j)}-\lambda_k | \leq
    |\lambda_k^{(j)}-\lambda_k(a_j,b) | +
    |\lambda_k(a_j,b)-\lambda_k|. \label{eq:mm24c} 
\ee
Theorem \ref{Theorem:FApprox2} gives
\be \lim_{j\rightarrow\infty}|\lambda_k(a_j,b)-\lambda_k| = 0
 \label{eq:mm25}
\ee
so it remains only to show that we can choose the $b_j$ 
and associated boundary conditions to get
\be \lim_{j\rightarrow\infty}|\lambda_k^{(j)}-\lambda_k(a_j,b) | = 0.
 \label{eq:mm26} 
\ee
We do this by noting that we are now dealing with problems
having just one singular lim-3 endpoint at $x=b$, the point
$x=a_j$ being regarded as regular and fixed. Provided the
appropriate hypotheses are satisfied, we shall be able to
apply Theorem \ref{Theorem:1}.

We denote by $\mu_k(a_j,b)$ the eigenvalues of the problems over
$(a_j,b)$ with Dirichlet conditions at $x=a_j$ and Friedrichs
conditions at the singular endpoint $x=b$. From Theorem \ref{Theorem:FApprox2}
we have not only (\ref{eq:mm25}) but also
\be \lim_{j\rightarrow\infty}|\mu_k(a_j,b)-\mu_k | = 0.
 \label{eq:mm27} 
\ee
Thus the inequality $\lambda_n < \mu_n$ translates to
$\lambda_n(a_j,b) < \mu_n(a_j,b)$ for all sufficiently large
$j$. This allows us to apply Theorem \ref{Theorem:1} to deduce that
there is a point $b_j < b$ and associated boundary conditions at
$x=b_j$, including the condition $[y,\psi](b_j)=0$, such that
\[ |\lambda_k^{(j)}-\lambda_k(a_j,b) | < \frac{1}{j}. \]
Combining this with (\ref{eq:mm25}) and (\ref{eq:mm24c}), our
proof is complete. \hfill $\Box$

\subsection{The lim-2, lim-4 case}
In this section we consider a problem posed over an interval $(a,b)$ in
which $x=a$ is a lim-2 endpoint and $x=b$ is lim-4. At $x=a$ we have
no boundary conditions. At $x=b$ we shall impose two boundary conditions
\be [y,u_3](b) = 0, \;\;\; [y,u_4](b) = 0 \label{eq:mm15} \ee
where, by recourse to the earlier results, we may assume that $u_3$ and
$u_4$ are, for some real $\lambda$, solutions of the differential equation
which are square integrable at $b$, linearly independent relative
to the minimal domain, and satisfy 
$[u_3,u_3] = [u_3,u_4] = [u_4,u_4] = 0$. As usual we
denote the resulting self-adjoint operator by $L$ and its eigenvalues
by $\lambda_k$.

We shall consider a sequence of approximating operators $L_j$ defined over
intervals $(a_j,b_j)$ where $a_j \rightarrow a$ and $b_j\rightarrow b$ as
$j\rightarrow \infty$. The boundary conditions (\ref{eq:mm15}) will now
be replaced by
\be y(a_j)=0=y'(a_j), \;\;\; [y,u_3](b_j) = 0 = [y,u_4](b_j). 
                                                    \label{eq:mm16} \ee
We shall denote the eigenvalues of these approximating problems by
$\lambda_k(a_j,b_j)$.

Following the results
which we presented in the case of just one lim-2 endpoint it is clear that
for spectral exactness we had better not have spectrum extending
to $-\infty$ with a gap. With a lim-4 endpoint we could, of course, have
discrete spectrum extending to $-\infty$. We shall assume that this is
not the case. Our main result is the following.

\begin{Theorem}
Suppose that $L$ is bounded below and possesses at least $n+1$ eigenvalues
$\lambda_0,\ldots,\lambda_n$ strictly below any essential spectrum. Suppose
also that the sequence $(a_j)_{j=1}^{\infty}$ converges monotonically to $a$.
Then we have the spectrally exact convergence
\be \lim_{j\rightarrow\infty}\lambda_k(a_j,b_j) = \lambda_k, \;\;\;
 k = 0,\ldots,n. \label{eq:mm17} \ee
\end{Theorem}
\noindent {\bf Proof}:\, To establish this result we compare the following
 eigenvalues:
\begin{itemize}
\item the eigenvalues $\lambda_k(a_j,b_j)$;
\item the eigenvalues $\lambda_k(a_i,b_j)$ for problems on intervals 
$(a_i,b_j)$;
\item the eigenvalues $\lambda_k(a,b_j)$ of problems with a lim-2 singular
 end at $x=a$ and a regular end at $x=b_j$, with boundary conditions 
$[y,u_3](b_j) = 0 = [y,u_4](b_j)$;
\item the eigenvalues $\lambda_k(a,b) := \lambda_k$ of $L$.
\end{itemize}
We start with the triangle inequality, which gives
\be |\lambda_k(a_j,b_j) - \lambda_k(a,b)| \leq 
 |\lambda_k(a_j,b_j)-\lambda_k(a,b_j)| + |\lambda_k(a,b_j)-\lambda_k(a,b)|
\label{eq:mm18} 
\ee
Next, we observe that for all $i < j$, we have
\be \lambda_k(a,b_j) \leq \lambda_k(a_j,b_j) \leq \lambda_k(a_i,b_j). 
\label{eq:mm19} \ee
This is because $a < a_j < a_i$ and the boundary conditions at the left
hand endpoints are always Dirichlet. This inequality yields
\be |\lambda_k(a_j,b_j) - \lambda_k(a,b_j)| \leq  
 |\lambda_k(a_i,b_j) - \lambda_k(a,b_j)| \label{eq:mm20} \ee
Substituting (\ref{eq:mm20}) into the right hand side of (\ref{eq:mm18}) yields
\be |\lambda_k(a_j,b_j) - \lambda_k(a,b) | \leq  
 |\lambda_k(a_i,b_j) - \lambda_k(a,b_j)|
 + |\lambda_k(a,b_j)-\lambda_k(a,b)| \;\; \forall j > i. \label{eq:mm21} \ee
If we now let $j$ tend to infinity on the right hand side of (\ref{eq:mm21})
then the second term, $|\lambda_k(a,b_j)-\lambda_k(a,b)|$, will tend to zero
by the results of Section \ref{section:lim4}, while the first term will tend
to $|\lambda_k(a_i,b) - \lambda_k(a,b)|$. Thus we have
\be \lim_{j\rightarrow\infty} |\lambda_k(a_j,b_j) - \lambda_k(a,b) | \leq 
 |\lambda_k(a_i,b) - \lambda_k(a,b)| \;\;\; \forall i. \label{eq:mm22b} \ee
To complete the proof one need only show that the right hand side of 
(\ref{eq:mm22b}) tends to zero as $i\rightarrow\infty$. This is an immediate
consequence of Theorem \ref{Theorem:FApprox2}. \hfill $\Box$

\subsection{The lim-3, lim-3 case}
Since the case of one lim-3 endpoint was so awkward it should come
as no surprise that the case of two lim-3 endpoints is the most
difficult to treat. We shall consider an eigenvalue problem 
\be \left\{ \begin{array}{l} 
            \ell y =  \lambda y, \;\;\; x \in (a,b), \\
             \left[ y,\phi\right](a) = 0, \\
             \left[ y,\psi\right](b) = 0, 
            \end{array} \right. \label{eq:mm40} 
\ee
which we denote by EP. We assume that the minimal operator
associated with this problem is bounded below,  that $\phi$ and $\psi$
are real, and that neither
$\phi$ nor $\psi$ defines a Friedrichs boundary condition (since
Friedrichs boundary conditions have been treated earlier).
We shall denote the eigenvalues of EP by $\lambda_k$. We also
require the problem
\be \left\{ \begin{array}{l} 
            \ell y =  \lambda y, \;\;\; x \in (a,b), \\
             \left[y,\phi\right](a) = 0, \\
             \mbox{Friedrichs BC at $b$.} 
            \end{array} \right. \label{eq:mm41} 
\ee
which we denote by $\mbox{EP}^{F}$ and whose eigenvalues we
denote by $\mu_k^F$, and the problem
\be \left\{ \begin{array}{l} 
            \ell y =  \lambda y, \;\;\; x \in (a,b), \\
             \mbox{Friedrichs BC at $a$,} \\
             \mbox{Friedrichs BC at $b$.} 
            \end{array} \right. \label{eq:mm42} 
\ee
which we denote by $\mbox{EP}^{FF}$ and whose eigenvalues we denote 
by $\mu_k^{FF}$. By variational methods it is easy to show that when
the relevant eigenvalues all exist and lie below any essential 
spectrum, then
\[ \lambda_k \leq \mu_{k}^{F} \leq \mu_k^{FF}. \]
In line with our assumptions in Section \ref{section:onelim3}, we
shall assume that for some positive integer $n$,
\be \lambda_n < \mu_n^F < \mu_n^{FF}. \label{eq:mm43} \ee
Our main result is that if, for some $n$, the inequalities
in (\ref{eq:mm43}) hold, then we can obtain a sequence of
approximating regular problems to give spectral exactness
for the first $n+1$ eigenvalues. 

\begin{Theorem}
Suppose that for some integer $n$, the eigenproblems EP,
$\mbox{EP}^F$ and $\mbox{EP}^{FF}$ all have at least $n+1$
eigenvalues strictly below the essential spectrum, and 
suppose that (\ref{eq:mm43}) holds. Then given any sequence
$(b_j)$ with $b_j\nearrow b$ as $j\nearrow\infty$,
we can construct a sequence $(a_j)$ with $a_j\searrow a$
as $j\nearrow\infty$, and a sequence of regular problems
$\mbox{EP}(a_j,b_j)$ on the intervals $[a_j,b_j]$ whose boundary conditions
include $[y,\phi](a_j)=0$ and $[y,\psi](b_j)=0$, such that
for $k=0,\ldots,n$, the $k$th eigenvalue $\lambda_k^{(j)}$
of $\mbox{EP}(a_j,b_j)$ satisfies
\[ \lim_{j\rightarrow\infty}\lambda_k^{(j)} = \lambda_k. \]
\end{Theorem}
\noindent {\bf Proof}:\, Since (\ref{eq:mm43}) holds, we
can fix a point $\lambda_*$ in $(\lambda_n,\mu_n^F)$ such that
$\lambda_* \neq \mu_k^F$ for all $k$.
Now given the sequence $b_j$ consider the sequence of
singular eigenproblems $\mbox{EP}^{F}(a,b_j)$ defined
by
\be \left\{ \begin{array}{l} 
            \ell y =  \lambda y, \;\;\; x \in (a,b), \\
             \left[y,\phi\right](a) = 0, \\
             y(b_j)=0=y'(b_j), 
            \end{array} \right. \label{eq:mm44} 
\ee
whose eigenvalues we denote $\mu_k^{F}(a,b_j)$. From Section
\ref{section:Friedrichs} we know that
\be \mu_k^{F}(a,b_j) \searrow \mu_k^F \;\;\; \mbox{as $j\nearrow\infty$},
 \;\;\; k=0,\ldots,n.
 \label{eq:mm45} \ee
Similarly, we can set up the eigenproblems $\mbox{EP}^{FF}(a,b_j)$ defined
by
\be \left\{ \begin{array}{l} 
            \ell y =  \lambda y, \;\;\; x \in (a,b), \\
             \mbox{Friedrichs BC at $a$,} \\
             y(b_j)=0=y'(b_j), 
            \end{array} \right. \label{eq:mm45b} 
\ee
whose eigenvalues we denote $\mu_k^{FF}(a,b_j)$. Again, Section
\ref{section:Friedrichs} tells us that
\be \mu_k^{FF}(a,b_j) \searrow \mu_k^{FF} \;\;\; 
\mbox{as $j\nearrow\infty$},
\;\;\; k=0,\ldots,n.
 \label{eq:mm46} \ee
Next, we examine $\mbox{EP}^F(a,b_j)$ for each $j$. This is a problem
with one lim-3 endpoint and one regular endpoint; we want to exploit
the results of Section \ref{section:onelim3} to approximate its 
eigenvalues. To this end we must check the hypothesis of Theorem
\ref{Theorem:1}: do we have $\mu_n^F(a,b_j)<\mu_n^{FF}(a,b_j)$?
Combining (\ref{eq:mm43}), (\ref{eq:mm45}) and (\ref{eq:mm46}), it
is clear that we do. Thus Theorem \ref{Theorem:1} can be applied:
we can choose an endpoint $\al_i>a$ and a set of boundary conditions
at $\al_i$ including the condition $[y,\phi](\al_i)=0$, such that the
eigenvalues, let us call them $\mu_k^F(\al_i,b_j)$, of the resulting
regular problem $\mbox{EP}^{F}(\al_i,b_j)$ approximate
the eigenvalues of $\mbox{EP}^{F}(a,b_j)$:
\be \lim_{i\rightarrow\infty} |\mu_k^F(\al_i,b_j) - \mu_k^F(a,b_j)| = 0,
\;\; k=0,\ldots,n. \label{eq:mm47} \ee
Choose $i = i(j)$ such that
\[  |\mu_k^F(\al_i,b_j) - \mu_k^F(a,b_j)| < \frac{1}{j}.  \]
Set $a_j = \al_{i(j)}$.  Then using (\ref{eq:mm45}) we have
\be \lim_{j\rightarrow\infty}\mu_k^F(a_j,b_j) = \mu_k^F,
\;\; k=0,\ldots,n. \label{eq:mm48}
\ee
The problems $\mbox{EP}^F(a_j,b_j)$ will have the form
\be \left\{ \begin{array}{l} 
            \ell y =  \lambda y, \;\;\; x \in (a,b), \\
             v(a_j) = W_L(a_j) u(a_j), \\
             y(b_j)=0=y'(b_j), 
            \end{array} \right. \label{eq:mm49} 
\ee
for some matrices $W_L(a_j)$ chosen so that the boundary condition
at $x=a_j$ is satisfied by $\phi$.

Given the problems $\mbox{EP}^F(a_j,b_j)$ we now aim to change the
boundary condition at $x=b_j$. With $W_L(a_j)$ fixed we attempt to
choose a matrix $W_R(b_j)$ such that the eigenproblems 
$\mbox{EP}(a_j,b_j)$ defined by
\be \left\{ \begin{array}{l} 
            \ell y =  \lambda y, \;\;\; x \in (a,b), \\
             v(a_j) = W_L(a_j) u(a_j), \\
             v(b_j) = W_R(b_j) u(b_j), 
            \end{array} \right. \label{eq:mm50} 
\ee
have precisely $n+1$ eigenvalues less than $\lambda_*$. We must do this
in such a way that the boundary condition at $x=b_j$ is satisfied by
$\psi$ in order to guarantee spectral inclusion. Choose $\epsilon > 0$
such that $[\lambda_*,\lambda_*+\epsilon] \subset [\lambda_*,\mu_n^{F})$.
For each $\lambda\in [\lambda_*,\lambda_*+\epsilon]$ we can integrate
the Hamiltonian form of the differential equation forward from $x=a_j$,
starting with initial conditions $U_L(a_j,\lambda) = I$ and $V_L(a_j,\lambda)
= W_L(a_j)$, to obtain $U_L(b_j,\lambda)$ and $V_L(b_j,\lambda)$.

We assert
that if $\epsilon$ is chosen sufficiently small then we shall have
$\det U_L(b_j,\lambda) \neq 0$ for all sufficiently large $j$ and for all 
$\lambda  \in [\lambda_*,\lambda_*+\epsilon]$.  For if this were not
true then we could extract a subsequence of the $b_j$ tending to $b$
and a sequence of values of $\lambda$ tending to $\lambda_*$ at which
we had $\det U_L = 0$. These values of $\lambda$ would be eigenvalues of
a subsequence of the problems $\mbox{EP}^{F}(a_j,b_j)$ and would
therefore have to converge to $\mu_k^{F}$ for some $k$, by (\ref{eq:mm48}).
This would mean that we had $\lambda_* = \mu_k^F $ for some $k$, which would
contradict the choice of $\lambda_*$.

If $u_{\psi}(b_j)\neq 0$ for sufficiently large $j$, then
Lemma \ref{lemma:sigma1} implies that we can choose the matrix
$W_R(b_j)$ so that the boundary condition $v(b_j)=W_R u(b_j)$ is
satisfied by $\psi$ and gives the eigenproblem $\mbox{EP}(a_j,b_j)$
precisely one more eigenvalue below $\lambda_*$ than
the eigenproblem $\mbox{EP}^F(a_j,b_j)$: in other words,
$\mbox{EP}(a_j,b_j)$ has the same number of eigenvalues below
$\lambda_*$ as $\mbox{EP}(a,b)$. Spectral inclusion (Lemma
\ref{lemma:mmsi}) now forces the convergence
\[ \lim_{j\rightarrow\infty}\lambda_k^{(j)} = \lambda_k, \;\;\;
 k = 0,\ldots, n. \]
If we have $u_{\psi}(b_j)=0$ on a subsequence of the $b_j$ then, as 
in the case of one lim-3 endpoint, we can show that the eigenvalues 
$\lambda_k^{(j)}$ must converge to $\mu_k^F$ for $k=0,\ldots, n$.
By spectral inclusion, together with the inequality $\lambda_k
\leq \mu_k^F$, this then implies that $\lambda_k = \mu_k^F$ for
$k=0,\ldots,n$, contradicting the assumption (\ref{eq:mm43}).
This completes the proof. \hfill $\Box$

\subsection{The lim-3, lim-4 case}
Suppose that $x=a$ is a lim-3 endpoint, $x=b$ is lim-4, and $L_{min}$ is
bounded below.  Let $L$ be a self-adjoint extension of $L_{min}$ with
boundary conditions $[y,\phi](a) = 0, \ \ [y,\psi_1](b) = 0 = [y,\psi_2](b)$.
We may assume that $\phi$ is real and $\psi_1,\ \psi_2$ are solutions of the
differential equation with real $\lam$.
We are interested in approximating the eigenvalues of $L$ below the essential
spectrum.  Let $L_0$ be the extension of $L_{min}$ with boundary conditions
$[y,\psi_1](b) = 0 = [y,\psi_2](b)$, and let $M$ be the Friedrichs extension
of $L_0$.  Suppose that $M$ has boundary conditions 
$[y,\th](a) = 0, \ \ [y,\psi_1](b) = 0 = [y,\psi_2](b)$.
Let $\lam_0 \le \lam_1 \le \lam_2 \le \cdots$ be the eigenvalues
of $L$ below the essential spectrum, and let $\mu_0 \le \mu_1 \le \mu_2 
\le \cdots$ be those of $M$.  Then $\lam_k \le \mu_k$.

\begin{Theorem}
Suppose that $\lam_n < \mu_n$ for some $n$.  Then we can construct a sequence
of regular operators $L_j$ on truncated intervals $[a_j,b_j]$, with eigenvalues
$\lam_k^{(j)}$, such that
\[ \lim_{j\rightarrow\infty}\lam_k^{(j)} = \lam_k, \;\;\; k = 0,1,\ldots,n. \]
\end{Theorem}
\vskip 5pt

\noindent {\bf Proof}:\,  Let $a < b_j < b, \ \ b_j \nearrow b \ {\rm as} \ 
j \rightarrow \infty$, and consider the following operators on the truncated
interval $(a,b_j]$:  $\hat{L}_j$ has boundary conditions 
$[y,\phi](a) = 0, \ \ [y,\psi_1](b_j) = 0 = [y,\psi_2](b_j)$;  $\hat{M}_j$ has
boundary conditions $[y,\th](a) = 0, \ \ [y,\psi_1](b_j) = 0 = [y,\psi_2](b_j)$.
$\hat{L}_j$ is a truncation of $L$, and $\hat{M}_j$ is a truncation of $M$.
By Corollary \ref{corollary:halftrunc}, the sequences $\hat{L}_j$ and 
$\hat{M}_j$
are spectrally exact for $L$ and $M$, respectively.  All of these operators have
the same essential spectrum.  Let $\hat{\lam}_k^{(j)}$ be the eigenvalues below
the essential spectrum for $\hat{L}_j$, and let $\hat{\mu}_k^{(j)}$ be those
for $\hat{M}_j$.  Thus
\be \lim_{j\rightarrow\infty}\hat{\lam}_k^{(j)} = \lam_k, \hm
     \lim_{j\rightarrow\infty}\hat{\mu}_k^{(j)} = \mu_k.  \label{eq:lg80}  \ee
Since $\lam_n < \mu_n$ it follows that for sufficiently large $j$
\be  \hat{\lam}_n^{(j)} < \hat{\mu}_n^{(j)}.  \label{eq:lg81}  \ee
Let $L_j^0$ be the operator on $(a,b_j]$ with boundary conditions
$[y,\psi_1](b_j) = 0 = [y,\psi_2](b_j)$; then $\hat{M}_j$ is the Friedrichs
extension of $L_j^0$.  Thus the hypothesis of Theorem \ref{Theorem:1} is
satisfied for $\hat{L}_j$.  Consequently, for a sequence $\al_i \searrow a$, 
we can construct regular operators $L_{i,j}$ on intervals $[\al_i,b_j]$, with
eigenvalues $\lam_k^{(i,j)}$, such that
\be  \lim_{i\rightarrow\infty}\lam_k^{(i,j)} = \hat{\lam}_k^{(j)}, 
                           \;\;\; k = 0,1,\ldots,n.     \label{eq:lg83}  \ee
For each $j$, choose $i = i(j)$ such that
\be   |\lam_k^{(i,j)} - \hat{\lam}_k^{(j)}| < \frac{1}{j}, 
                       \;\;\; k = 0,1,\ldots,n.     \label{eq:lg84}  \ee
Let $a_j = \al_{i(j)}, \ L_j = L_{i(j),j}, \ 
\lam_k^{(j)} = \lam_k^{(i(j),j)}.$  Then $L_j$ is a regular operator on
the interval $[a_j,b_j]$ with eigenvalues $\lam_k^{(j)}$, and
$\ds \lim_{j\rightarrow\infty}\lam_k^{(j)} = \lam_k$, 
for $k = 0, 1, 2, \ldots, n.$  \hfill $\Box$

\end{document}